\documentclass[11pt,dvipsnames]{article}
\usepackage[T1]{fontenc}

\usepackage[a4paper,margin=2.8cm,top=2.2cm,bottom=2cm]{geometry}

\usepackage{amsmath}
\usepackage{array}
\usepackage{booktabs}
\usepackage{caption}
\usepackage{fancyvrb}
\usepackage{graphicx}
\usepackage{hyperref}
\usepackage{url}
\hypersetup{breaklinks=true, colorlinks,
  urlcolor=RoyalBlue, citecolor=PineGreen, linkcolor=RoyalBlue}
\usepackage{libertine}
\usepackage{longtable}
\usepackage{newtxmath}
\usepackage{tabularx}
\usepackage{zi4}

\usepackage[english]{babel}
\usepackage{enumitem}
\usepackage{graphicx,subfigure,bm}
\usepackage[ruled,algo2e,linesnumbered]{algorithm2e}
\DontPrintSemicolon
\usepackage{myshortcuts}
\usepackage[utf8]{inputenc}

\usepackage{fancyvrb}

\usepackage{jlcode}

\pgfplotsset{compat=1.5.1}

\newtheorem{theorem}{Theorem}
\newtheorem{remark}[theorem]{Remark}
\newtheorem{example}[theorem]{Example}

\usepackage{tikz,pgfplots,pgfplotstable}
\usepackage{pgf}
\usetikzlibrary{decorations.pathreplacing}

\iftrue 
  \usepackage{tikz}
  \usepackage{pgfplots}
  \pgfplotsset{
    compat=newest,
    tick label style={font=\scriptsize},
    label style={font=\scriptsize},
    legend style={font=\scriptsize}
  }
  \iffalse 
     \usepgfplotslibrary{external} 
  \else
     \renewcommand{\tikzsetnextfilename}[1]{}

  \fi
\else
  \usepgfplotslibrary{external} 
  \tikzsetexternalprefix{gfx_tmp/}
  
\fi

\pgfplotscreateplotcyclelist{bwd-list}{%
  every mark/.append style={solid}, mark=none, solid, Cyan \\
  every mark/.append style={solid}, mark=none, dash pattern=on 0.8 off 0.8, Orange\\
  every mark/.append style={solid}, mark=x, mark size=4pt, solid, Black\\
}

\pgfplotscreateplotcyclelist{bwd-list-2}{%
  every mark/.append style={solid}, mark=none, solid, Red \\
  every mark/.append style={solid}, mark=none, densely dashed, Plum \\
  every mark/.append style={solid}, mark=none, dashdotted, Green \\
  every mark/.append style={solid}, mark=none, dash pattern=on 0.8 off 0.8 on
  0.8 off 0.8 on 3.2 off 0.8, Cyan \\
  every mark/.append style={solid}, mark=none, dash pattern=on 0.8 off 0.8, Orange \\
}

\pgfplotscreateplotcyclelist{bwd-list-3}{%
  every mark/.append style={solid}, mark=none, solid, Red \\
  every mark/.append style={solid}, mark=none, densely dashed, NavyBlue \\
  every mark/.append style={solid}, mark=none, dashdotted, Green \\
  every mark/.append style={solid}, mark=none, densely dotted, Red \\
  every mark/.append style={solid}, mark=none, densely dotted, NavyBlue \\
}

\newcommand{\thetitle}{Computational graphs for matrix functions}
\newcommand{\thetitlenote}{%
  Version of \today.
  The work of the second author was supported by the Wenner-Gren Foundations through grant UPD2019-0067.
}
\title{\huge\bfseries\thetitle\footnote{\thetitlenote}}

\newcommand{\theauthori}{Elias Jarlebring}
\newcommand{\theaffiliationi}{%
  \department{Department of Mathematics}
  \institution{KTH Royal Institute of Technology,
    SeRC Swedish e-science research center}
  \streetaddress{Lindstedtsvägen 25}
  \city{Stockholm}
  \postcode{SE-100 44}
  \country{Sweden}}
\newcommand{\theemaili}{eliasj@kth.se}
\newcommand{\theauthorii}{Massimiliano Fasi}
\newcommand{\theaffiliationii}{%
  \department{School of Science and Technology}
  \institution{Örebro University}
  \city{Örebro}
  \postcode{SE-701 82}
  \country{Sweden}}
\newcommand{\theemailii}{massimiliano.fasi@oru.se}

\newcommand{\theauthoriii}{Emil Ringh}
\newcommand{\theaffiliationiii}{%
  \department{Department of Mathematics}
  \institution{KTH Royal Institute of Technology,
    SeRC Swedish e-science research center}
  \streetaddress{Lindstedtsvägen 25}
  \city{Stockholm}
  \postcode{SE-100 44}
  \country{Sweden}}
\newcommand{\theemailiii}{eringh@kth.se}

\newcommand{\department}[1]{#1,}
\newcommand{\institution}[1]{#1,}
\newcommand{\streetaddress}[1]{#1,}
\newcommand{\city}[1]{#1}
\newcommand{\postcode}[1]{#1,}
\newcommand{\country}[1]{#1}
\author{\theauthori\thanks{\theaffiliationi, \href{mailto:\theemaili}{\theemaili}.} \and%
  \theauthorii\thanks{\theaffiliationii, \href{mailto:\theemailii}{\theemailii}.} \and %
  \theauthoriii\thanks{\theaffiliationiii, \href{mailto:\theemailiii}{\theemailiii}.}}%
\date{}

\newcommand{\theabstract}{%
  Many numerical methods for evaluating matrix functions can be naturally viewed as computational graphs. Rephrasing these methods as directed acyclic graphs (DAGs) is a particularly effective approach to study existing techniques, improve them, and eventually derive new ones.
  The accuracy of these matrix techniques can be characterized by the accuracy of their scalar counterparts, thus designing algorithms for matrix functions can be regarded as a scalar-valued optimization problem. The derivatives needed during the optimization can be calculated automatically by exploiting the structure of the DAG, in a fashion analogous to backpropagation.
  This paper describes \texttt{GraphMatFun.jl}, a Julia package that offers the means to generate and manipulate computational graphs, optimize their coefficients, and generate Julia, MATLAB, and C code to evaluate them efficiently at a matrix argument. The software also provides tools to estimate the accuracy of a graph-based algorithm and thus obtain numerically reliable methods. For the exponential, for example, using a particular form (degree-optimal) of polynomials produces implementations that in many cases are cheaper, in terms of computational cost, than the Padé-based techniques typically used in mathematical software.
  The optimized graphs and the corresponding generated code are available online.
}

\newcommand{\thekeywords}{%
  Polynomials of matrices,
  functions of matrices,
  computational graphs
}

\begin{document}
\maketitle

\begin{abstract}
  \noindent\theabstract
  \bigskip

\end{abstract}

\section{Introduction}

The scalar function $f : \Omega \subset \CC \to \CC$ can be extended to square
matrices in a customary fashion.
Formally, $f(A)$ can be defined in a number of ways: one can use, for example, the
Jordan decomposition, a Taylor series expansion, Cauchy integrals, or Hermite interpolation~\cite[Chapter 1]{high08}.
For a matrix $A$ whose spectrum lies in a region where $f$ is analytic, these four  definitions are equivalent. Matrix functions play an important role in numerical linear algebra as well as in matrix theory, and the existing literature ranges from theoretical results to computation and applications~\cite{high08}.

Numerical algorithms for evaluating a functions $f$ at an $\nbyn$ matrix $A$ are based on a variety
of different approaches, but for a large class of methods the approximation is
formed by relying on three basic operations:
\begin{enumerate}[label=O\arabic*.,ref=O\arabic*]
\item\label{M1} linear combination of matrices,
\item\label{M2} matrix multiplication, and
\item\label{M3} solution of a linear systems with $n$ right-hand sides.
\end{enumerate}
Our work focuses on algorithms constructed by applying \ref{M1}--\ref{M3}.
This approach does not include methods based on the Schur decomposition, such
as the algorithm of Bj\"orck and Hammarling for the square root~\cite{bjha83,
  high87s} or the Schur--Parlett algorithm for analytic matrix
functions~\cite{parl76}, nor does it include
methods for computing the action of a function
on a vector \cite{alhi11} or for computing individual matrix elements \cite{bsh21}.
Restricting oneself to~\mbox{\ref{M1}--\ref{M3}} yields matrix functions that
are particularly easy to work with: loosely speaking, if the expression for
$f(z)$ features only the scalar counterparts of~\ref{M1}--\ref{M3} and the
function is defined on the spectrum of~$A$~\cite[Definition~1.1]{high08}, then
a formula for $f(A)$ can be obtained by simply replacing all occurrences of $z$
in the formula for $f(z)$ with $A$~\cite[Chapter~9]{gova13}.

By interpreting functions constructed in this way as computational graphs where every node is a matrix, we can derive new methods and improve the
accuracy and performance of state-of-the-art algorithms for the evaluation of matrix functions.
The contribution of this work is twofold.
\begin{enumerate}
\item We present a new software package, \texttt{GraphMatFun.jl}, designed to
  work with computational graphs featuring only~\ref{M1}--\ref{M3}.
  The package includes tools to create and manipulate such graphs, modify their
  topology, optimize the matrix functions they represent, and evaluate them at
  scalar, vector, and matrix arguments.
  A framework to generate efficient Julia, MATLAB, and C code that implements
  these functions is also provided.
  The software is written in Julia
  ~\cite{beks17}, is released under the terms of the MIT license, and is
  freely available on GitHub.\footnote{\url{https://github.com/matrixfunctions/GraphMatFun.jl}}
  A detailed documentation is available via the built-in Julia documentation system and as part of the online user reference manual.\footnote{\url{https://matrixfunctions.github.io/GraphMatFun.jl/dev}}
  Pull
  requests are welcome, and user contribution is facilitated by the availability
  of a complete set of unit tests and of a continuous integration infrastructure.

\item We provide a set of graphs that were obtained by using our optimization
  strategy on algorithms available in the literature.
  In many cases of practical interest, the resulting algorithms are more
  efficient and accurate than those they originated from, and in fact of any
  existing alternative.
  These specific graphs represent a scientific contribution in their own right and are available on GitHub in a separate data repository.\footnote{\url{https://github.com/matrixfunctions/GraphMatFunData}}
\end{enumerate}%

The features of the package are presented by means of examples accompanied by
extended code snippets.
\begin{itemize}
\item The facilities for graph generation and basic node handling are illustrated in Examples~\ref{exmp:firstexample} and~\ref{exmp:addfunctions}.
\item The data structure that underpins the graph and the numerical evaluation of the function underlying a graph are illustrated in Examples~\ref{exmp:topoorder} and~\ref{exmp:firstexample}.
\item The graph manipulation and automatic optimization techniques are illustrated in Example~\ref{exmp:compress}.
\item The degree-optimal form for evaluating polynomial with a reduced number of non-scalar multiplications is
  described in Section~\ref{sec:degopt} and illustrated in Example~\ref{exmp:degopt}.
\item The features associated with the use of optimization to design graph algorithms are outlined in Section~\ref{sec:opt} and illustrated in Examples~\ref{exmp:Jac} and~\ref{exmp:optim}.
\item The generation of code that evaluates efficiently the function underlying a graph is discussed in Section~\ref{sec:codegen} and illustrated in Example~\ref{exmp:codegen}.
\item The routines for saving and loading graph objects in a language-independent format are introduced in Section~\ref{sect:cgr} and illustrated in Example~\ref{exmp:cgr}.
\item The backward error analysis of methods for approximating the matrix exponential is given in Section~\ref{sec:bwerrexp} and illustrated in Example~\ref{exmp:bwerr}.
\item A simple strategy to compute running error bounds on the round-off error occurring during the evaluation of a graph are given in Section~\ref{sec:roundoff} and illustrated in Example~\ref{exmp:roundoff}.
\end{itemize}

While developing our package, we realized that some of the Julia \texttt{Base} routines
for evaluating scalar polynomials and for computing the exponential of a matrix were not always optimal in terms of performance.
We addressed these issues in pull requests that have now been merged into the master branch of the Julia GitHub
repository.\footnote{\url{https://github.com/JuliaLang/julia/pulls?q=author:jarlebring}} Code generated by this package has been included in the Julia package \verb#SciML/ExponentialUtilities.jl#.

Section~\ref{sec:algorithms-as-graphs} introduces our graph framework and shows how our package can be used to generate and manipulate graphs that represent matrix algorithms.
In Section~\ref{sec:degopt} we discuss the polynomials in degree-optimal form, which are the most general polynomials of degree $2^m$ that can be evaluated with $m$ multiplications.
Degree-optimal polynomials can be represented in a natural way using our graph framework.
Optimizing the coefficients of a graph so that the underlying algorithm approximates a function that cannot be expressed as a computational graph is the subject of Section~\ref{sec:opt}.
Section~\ref{sec:eff-rep} is devoted to the description of two important features of the package: the generation of Julia, MATLAB, and C code for evaluating a computational graph efficiently, and the custom file format for storing computational graphs.
Using a graph to approximate a function numerically is prone to truncation as well as round-off errors; these are discussed in Section~\ref{sec:erroranalysis}, where we obtain a priori bounds on the former and a posteriori running bounds on the latter.
In Section~\ref{sec:use-cases} we explain how the features described in the previous sections can be integrated into a workflow for enhancing state-of-the-art algorithms for matrix functions, and present two instances in which such improvement can be encountered in practice.
Finally, we summarize our contribution and outline possible directions for future investigation in Section~\ref{sec:conclusions}.

\section{Algorithms as graphs}
\label{sec:algorithms-as-graphs}

\subsection{Graph notation and operations}
We regard matrix computations as graphs,
in particular as directed acyclic graphs (DAGs), in
which every node corresponds to a matrix.
The leaves of the graph coincide with the input nodes,
and the corresponding matrices
are assumed to be known in advance.
Every non-input node is constructed by
performing an operation on two other matrices which can be
identified by following the two incoming edges.
Hence, input nodes have only outgoing edges, whereas every
non-input node has exactly two incoming edges and
zero or more outgoing edges. The output nodes represent the
output (or outputs) of the computation.

Although many operations could potentially be considered in such a framework,
in this work and in our software we have restricted our attention to the following:
\begin{itemize}
\item $+$, which denotes the linear combination of two matrices, $C=\alpha A+\beta B$;
\item $*$, which denotes the product of two matrices, $C=AB$; and
\item $\backslash$, which denotes the solution of a linear systems with multiple right-hand sides, $C=A^{-1}B$.
\end{itemize}
With these three building blocks one can construct many algorithms of practical
interest for evaluating matrix functions. Moreover, having only
three operations all acting on exactly
two matrices has several implementation advantages.

Let $V$ be the set of nodes. For every node $i\in V$ we define:
\begin{itemize}
\item $p_i\in \{+,*,\backslash\}$, the operation required to evaluate node $i$;
\item $e_i^{(1)}\in V$ and $e_i^{(2)}\in V$, the two nodes representing the matrices needed compute node $i$; and
\item $c_i^{(1)}$ and $c_i^{(2)}$, the coefficients of the linear combination, defined only when $p_i=+$.
\end{itemize}
In order to evaluate the function underlying the computational graph, we need to explicitly construct the matrices each node represents.
To this end, we also define:
\begin{itemize}
\item $Z_i$, the content (matrix or scalar) of node $i$.
\end{itemize}
We will sometimes use the functional notation $Z_i(\cdot)$ to stress that $Z_i$ depends on the content of the input nodes.

\begin{figure}[t]
  \begin{center}
    \subfigure[Denman--Beavers iteration for the square root \cite{debe76}.]{\includegraphics[width=0.47\textwidth]{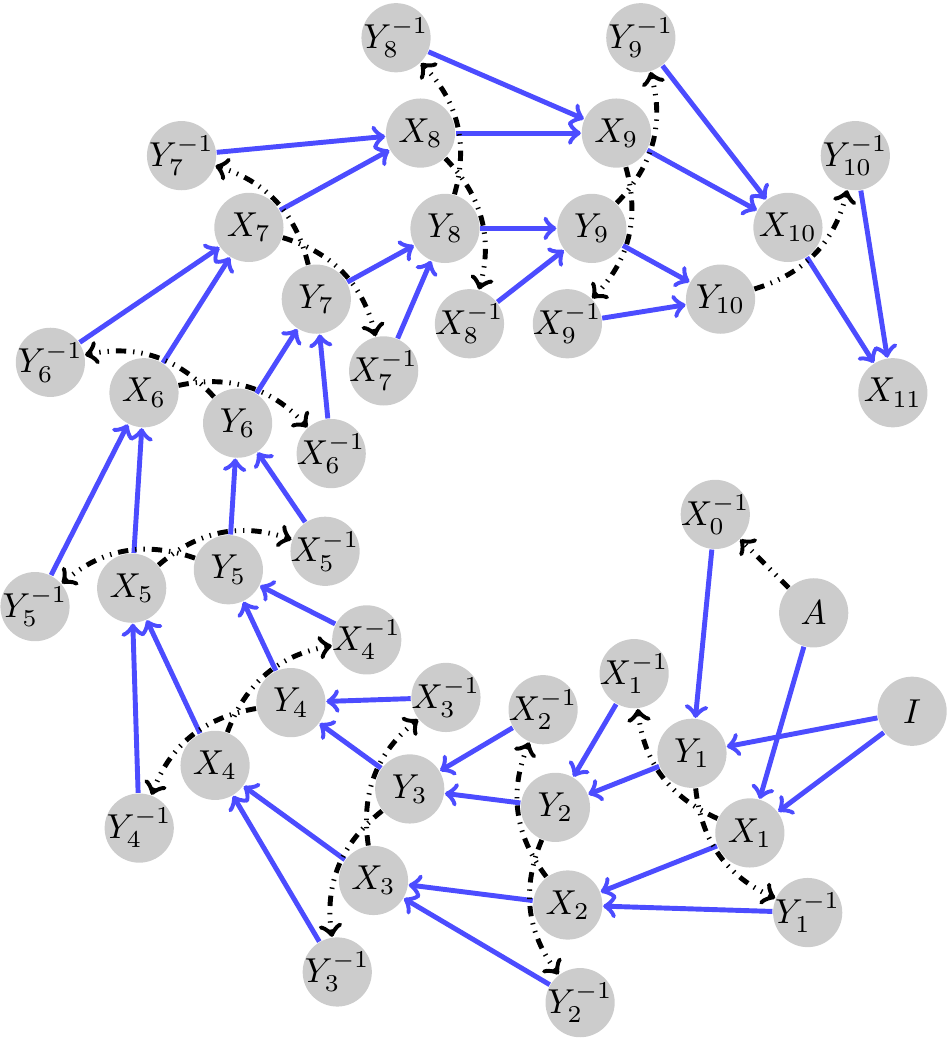}\label{fig:db}}%
    \qquad%
    \subfigure[Scaling and squaring for the matrix exponential~\cite{high05e}.]{\includegraphics[width=0.47\textwidth]{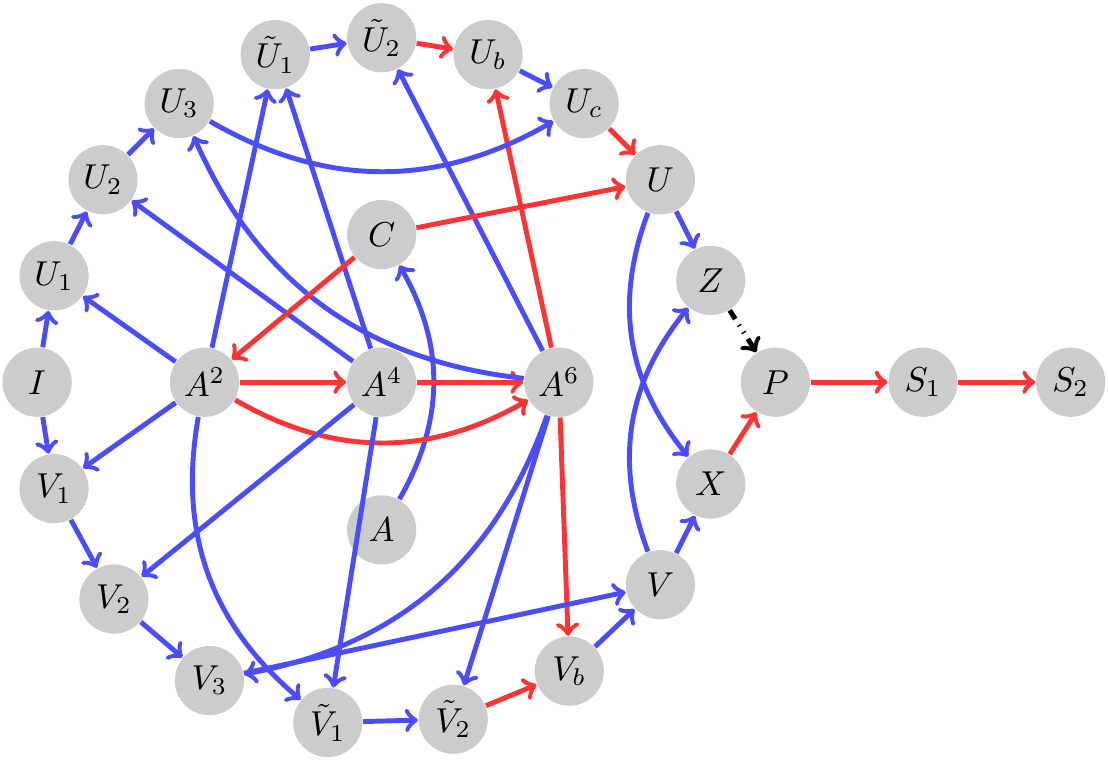}}
    \caption{Examples of complex graphs. Blue lines represent linear combinations, red lines represent matrix multiplications, and black dash-dotted lines represent system solves. Nodes with only one visible parent are either the result of a squaring ($A^2$, $A^4$, $S_1$, $S_2$), or a suppressed trivial operation involving the identity matrix: a matrix inversion, represented as  $C=A^{-1}I$, or a multiplication by a scalar, stored  as $C = \alpha A + 0 I$.
      \label{fig:advancedgraphs}
    }
  \end{center}
\end{figure}

For notational convenience, here we denote the nodes by integers and write $V=\{0,\ldots,N\}$. In the
software the set of nodes is implemented as an (unordered) dictionary lookup table for efficiency.
The two input nodes are the identity matrix and the $A$ matrix, which are labelled $Z_0=I$ and $Z_1=A$ in the matrix case, and  $Z_0=1$ and $Z_1=x$ in the scalar case.

The values of the other nodes are given by
\newsavebox{\casescommonpart}
\begin{subequations}\label{eq:fundop}
  \begin{align}
    \sbox{\casescommonpart}{$\displaystyle Z_i(x) =\left\{%
    \begin{array}{@{}l@{}}%
      \vphantom{c_i^{(1)} Z_\leftp(x) + c_i^{(2)} Z_\rightp(x)}\\[2pt]%
      \vphantom{Z_\leftp(x) Z_\rightp(x)}\\[2pt]%
      \vphantom{Z_\leftp(x)^{-1}Z_\rightp(x)}\\[2pt]%
    \end{array}%
    \right.\kern-\nulldelimiterspace$}
    \raisebox{-.66\ht\casescommonpart}[0pt][0pt]{\usebox{\casescommonpart}}%
    \parbox{3.375cm}{$c_i^{(1)} Z_\leftp(x) + c_i^{(2)} Z_\rightp(x),$} \qquad p_i=+,\label{eq:lincomb}\\
    \parbox{3.375cm}{$Z_\leftp(x) Z_\rightp(x),$}\qquad p_i=*, \label{eq:mult}\\
    \parbox{3.375cm}{$Z_\leftp(x)^{-1}Z_\rightp(x),$}\qquad p_i=\backslash, \label{eq:ldiv}
  \end{align}
\end{subequations}
where $\leftp=e^{(1)}_i$ and $\rightp=e^{(2)}_i$.

The set of (directed) edges of the DAG is
\[
\Bigl\{
\bigl(e_2^{(1)},2\bigr),\bigl(e_2^{(2)},2\bigr),
\bigl(e_3^{(1)},3\bigr),\bigl(e_3^{(2)},3\bigr),
\ldots,
\bigl(e_N^{(1)},N\bigr),\bigl(e_N^{(2)},N\bigr)\Bigr\}.
\]
These edges form a parent list. Therefore,
if we assume that the graph has no directed cycles,
by traversing it a in a depth-first fashion
we can find a topological ordering based on the parent-child
relationship.
A topological ordering is a linear arrangement of the nodes in which every parent
follows both of its children.
The topological ordering is not unique, and the maximum number of nodes that
are needed simultaneously at any stage is called the width of the
DAG. Computing the width and the corresponding topological
ordering is a difficult problem~\cite{bodl93,bgt98},
which arises for example, in compiler optimization.
In our setting, topological sorting will be of paramount importance,
since a good ordering of the nodes decreases
the amount of memory required to evaluate the graph at a matrix argument.

\begin{example}[Graph evaluation and topological ordering]\label{exmp:topoorder}
  The graph corresponding to the first four steps of the Denman--Beavers iteration~\cite{debe76} can be constructed with the command:
  \begin{lstlisting}
julia> using GraphMatFun
julia> (graph, _) = graph_denman_beavers(4); # 4 iterations
  \end{lstlisting}
  and can be evaluated at a matrix argument using \verb#eval_graph#:
  \begin{lstlisting}
julia> A = [0.5 0.2; 0.3 0.5];
julia> eval_graph(graph, A)
2×2 Matrix{Float64}:
 0.684065  0.146185
 0.219277  0.684065
julia> sqrt(A)  # Compare with the Julia's matrix square root
2×2 Matrix{Float64}:
 0.684065  0.146185
 0.219277  0.684065
\end{lstlisting}
  The graph is depicted in Figure~\ref{fig:db}.
  This evaluation relies on the ordering of the nodes computed by \verb#get_topo_order#, which can be examined by calling the latter directly.
  \begin{lstlisting}
julia> get_topo_order(graph)[1]
17-element Vector{Symbol}:
 :Xinv0
 :Y1
 :X1
 :Yinv1
 :X2
 ...
 :Xinv3
 :Y4
 :Yinv4
 :X5
  \end{lstlisting}
\end{example}

Figure~\ref{fig:advancedgraphs} shows the graphs produced by two algorithms of practical interest: the Denman--Beavers iteration for the square root~\cite{debe76} and the scaling and squaring algorithm for the matrix exponential based on a degree-13 Padé approximant~\cite{high05e}.

\subsection{Basic package features}
The package provides a variety of functions for creating and manipulating
graphs.
The routines to generate directly
graphs corresponding to several well-established algorithms
are listed in Table~\ref{tbl:generators}. Table~\ref{tbl:addfunctions} and
Table~\ref{tbl:graphfunctions} list functions that can be used to generate new
graphs and manipulate existing ones, respectively.

\begin{table}[p]
  \begin{center}
    \caption{Graph generators. For each routine we report the name, the function being evaluated/approximated, and a reference for the implementation. The degree-optimal form is discussed in Section~\ref{sec:degopt}. \label{tbl:generators}}
    \setlength{\tabcolsep}{3.338pt}
    \vskip -3pt
    \begin{tabular}{@{}lcl@{}}
      \toprule
      Routine name & $f(x)$  & Implemented algorithm \\
      \midrule
      \verb#graph_monomial#    & $p(x)$& Standard monomial evaluation \\
      \verb#graph_monomial_degopt#    & $p(x)$& Above embedded into degree-optimal form \\
      \verb#graph_horner#      & $p(x)$& Horner scheme  evaluation \\
      \verb#graph_horner_degopt#      & $p(x)$& Above embedded into degree-optimal form \\
      \verb#graph_ps#      & $p(x)$&  Paterson--Stockmeyer \cite{past73} \\
      \verb#graph_ps_degopt#      & $p(x)$&  Above embedded into degree-optimal form \\
      \verb#graph_sastre_poly#     & $p(x)$ &  Efficient evaluation of polynomial of degree 8 \cite{sast18}\\
      \verb#graph_sastre_yks_degopt#     & $p(x)$ &  Transforms $y_{ks}$-form \cite{sast18} to degree-optimal form \\
      \verb#graph_rational#     & $r(x)$ &  Rational with poly eval with any of above \\
      \verb#graph_denman_beavers#     & $\sqrt{x}$ & Denman--Beavers iteration \cite{debe76} \\
      \verb#graph_newton_schulz#     & $x^{-1}$ &  Newton--Schulz iteration \cite{schu33} \\
      \verb#graph_newton_schulz_degopt#     & $x^{-1}$ &  Above embedded into degree-optimal form \\
      \verb#graph_exp_native_jl#     & $e^x$ &  Scaling and squaring with Padé approximant \cite{high05e}\\
      \verb#graph_exp_native_jl_degopt#     & $e^x$ &  Above embedded into degree-optimal form\\
      \verb#graph_bbc_exp#     & $e^x$ &   Scaling and squaring with Taylor approximant \cite{bbc19}\\
      \verb#graph_bbcs_cheb_exp#     & $e^x$ &  Chebyshev approximant for skew-Hermitian matrices \cite{bbcs21}\\
      \verb#graph_sastre_exp#     & $e^x$ &  Efficient evaluation of Taylor approximant \cite{sast18}\\
      \verb#graph_sid_exp#     & $e^x$ &  Scaling and squaring with Taylor approximant \cite{sast18,sid19}\\
      \bottomrule
    \end{tabular}
  \end{center}
\end{table}

\begin{table}[p]
  \begin{center}
    \caption{Fundamental node manipulation functions. The \texttt{add\_sum!} function is a convenience function implemented by multiple calls to \texttt{add\_lincomb!}.
      \label{tbl:addfunctions}}
    \setlength{\tabcolsep}{11.25pt}
    \vskip -3pt
    \begin{tabular}{@{}ll@{}}
      \toprule
      Routine name & Operation \\
      \midrule
      \verb#add_lincomb!# & Adds a new node $Z_j=c_j^{(1)} Z_\leftp+c_j^{(2)} Z_\rightp$ \\
      \verb#add_mult!# & Adds a new node $Z_j=Z_\leftp Z_\rightp$ \\
      \verb#add_ldiv!# & Adds a new node $Z_j=Z_\leftp^{-1}Z_\rightp$ \\
      \verb#add_sum!# & Adds a new node $Z_j=\alpha_1Z_{m_1}+\cdots+\alpha_n Z_{m_n}$  \\
      \verb#del_node!# & Removes an existing node \\
      \verb#rename_node!# & Changes the label (\verb#Symbol#) of an existing node in all graph data structures\\
      \verb#add_output!# & Marks a node as output node \\
      \verb#clear_outputs!# & Clears the list of output nodes  \\
      \bottomrule
    \end{tabular}
  \end{center}
\end{table}

\begin{table}[p]
  \begin{center}
    \caption{Graph manipulation and evaluation functions.
      \label{tbl:graphfunctions}}
    \setlength{\tabcolsep}{8.55pt}
    \vskip -3pt
    \begin{tabular}{@{}ll@{}}
      \toprule
      Routine name & Operation \\
      \midrule
      \verb#Compgraph(T)# & Creates an empty graph with coefficients of type \verb#T#\\
      \verb#Compgraph(T,orggraph)# & Converts the graph \verb#orggraph# so that all its coefficients have type \verb#T# \\
      \verb#eval_graph# & Evaluates the graph at  a given matrix, scalar, or other object\\
      \verb#eval_jac# & Evaluates the Jacobian at a point (element-wise for a vector)\\
      \verb#compress_graph!# & Removes unnecessary operations from the graph \\
      \verb#merge_graph# & Combines two graphs into one and renames the nodes accordingly\\
      \bottomrule
    \end{tabular}
  \end{center}
\end{table}

\begin{figure}[t]
  \begin{center}%
\begin{Verbatim}[frame=single]
struct Compgraph{T}
    operations::Dict{Symbol, Symbol}
    parents::Dict{Symbol, Tuple{Symbol, Symbol}}
    coeffs::Dict{Symbol, Tuple{T, T}}
    outputs::Vector{Symbol}
end
\end{Verbatim}
    \caption{The fundamental data structure in the package
      is the \texttt{Compgraph}, which represents a graph
      and its associated operations. The node identifiers
      are \texttt{Symbol}s. As all operations involve
      two other nodes, the parent list is a dictionary (hash table)
      that associates a 2-element \texttt{Tuple} of parents to each
      non-input node.
    }\label{fig:datastructure}
  \end{center}
\end{figure}

\begin{example}[Graph]\label{exmp:firstexample}
  From a user's perspective, the easiest way to
  generate a graph is by using one of the graph
  generators in Table~\ref{tbl:generators}.
  The code corresponding to the graph can then be
  executed by calling \verb#eval_graph#.
\begin{lstlisting}
julia> (g_mono, _) = graph_monomial([1.0, 0.0, 3.0]); # Create graph for 1+3x^2
julia> x = 0.1;
julia> eval_graph(g_mono, x)
1.03
julia> 1 + 3x^2
1.03
julia> A = [3.0 4.0; 5.0 6.0];
julia> P = eval_graph(g_mono, A) # can be used for scalars or matrices
2x2 Matrix{Float64}:
  88.0  108.0
 135.0  169.0
\end{lstlisting}
The graph is a \verb#struct# with four fields (the definition of the data structure is given in Figure~\ref{fig:datastructure}):
\begin{lstlisting}
julia> g_mono.parents # Every node has two parents
Dict{Symbol, Tuple{Symbol, Symbol}} with 3 entries:
  :P3 => (:P2, :A2)
  :A2 => (:A, :A)
  :P2 => (:I, :A)
julia> g_mono.operations # Every node has an associated operation
Dict{Symbol, Symbol} with 3 entries:
  :P3 => :lincomb
  :A2 => :mult
  :P2 => :lincomb
julia> g_mono.coeffs # The linear combination operations have coefficients
Dict{Symbol, Tuple{Float64, Float64}} with 2 entries:
  :P3 => (1.0, 3.0)
  :P2 => (1.0, 0.0)
julia> get_topo_order(g_mono)[1] # A computable order
3-element Vector{Symbol}:
 :A2
 :P2
 :P3
\end{lstlisting}
\end{example}
\begin{example}[Compression]\label{exmp:compress}
  Node manipulation may lead to graphs whose topology can
  be optimized, for example
  by removing operations with zero coefficients.
  In Example~\ref{exmp:firstexample},
  for instance,
  the second coefficients of node \verb#:P2# is zero,
  since the term $x$ does not appear in the polynomial $1 + 3x^2$
  represented by \verb#g_mono#.
  We can therefore compress the graph and obtain a more efficient
  algorithm.
  We can achieve this by using \verb#compress_graph!#, a function that
removes 1) dangling nodes, i.e., non-output nodes that are not used anywhere in the graph; 2) trivial nodes, i.e., multiplication by the identity matrix or linear systems whose coefficient is the identity matrix; 3) redundant nodes, i.e., nodes that repeat computation already present in the graph; and 4) pass-through nodes, i.e., linear combinations whose coefficients include zero.
  \begin{lstlisting}
julia> using BenchmarkTools;
julia> A=randn(100, 100);
julia> @btime eval_graph($g_mono, $A);
  172.168 μs (191 allocations: 330.75 KiB)
julia> compress_graph!(g_mono);
julia> g_mono.parents # The node :P2 has been removed
Dict{Symbol, Tuple{Symbol, Symbol}} with 2 entries:
  :P3 => (:I, :A2)
  :A2 => (:A, :A)
julia> @btime eval_graph($g_mono, $A); # It runs faster with fewer allocations
  133.430 μs (123 allocations: 246.20 KiB)
  \end{lstlisting}
\end{example}

\begin{example}[Manipulating graph topology]\label{exmp:addfunctions}
  The following illustrates how a graph can be manipulated
  by adding nodes with the \verb#add_<operation># functions
  in Table~\ref{tbl:addfunctions}.
  \begin{lstlisting}
julia> g_mono.outputs # This is the output node
1-element Vector{Symbol}:
:P3
julia> add_mult!(g_mono, :PX, :P3, :P3); # Let's square the output node
julia> clear_outputs!(g_mono);
julia> add_output!(g_mono, :PX); # And set a new output node
julia> eval_graph(g_mono, x);
1.0609
julia> (1 + 3x^2)^2
1.0609
\end{lstlisting}
\end{example}

\section{Degree-optimal polynomials}\label{sec:degopt}

If $p$ is a polynomial, then evaluating $p(A)$ requires only matrix multiplications and linear combinations, and since the computational cost of the former asymptotically dominates that of the latter, it is crucial to use polynomial evaluation schemes that require as few matrix multiplications as possible.
The Paterson--Stockmeyer method~\cite{past73} requires fewer matrix products than the naive scheme in which the powers of $A$ are computed explicitly, and is currently the most efficient method for evaluating $p(A)$ using only the original monomial coefficients of $p$~\cite{fasi19}.
If the use of new coefficients derived from those of $p$ is allowed, however, then even more economical schemes, such as~\cite[Algorithm~C]{past73}, are possible.
Fuelled by the seminal work of Sastre~\cite{sast18}, the topic has received renewed attention in recent years, and algorithms that combine this scheme with the scaling and squaring technique have been developed for the exponential~\cite{bbc19,sid19}, the sine and cosine~\cite{siap19,sbbc21}, and the hyperbolic tangent~\cite{iasd21}.

Polynomials of degree $2^m$ can be obtained with $m$ multiplications by repeated squaring.
This is optimal, in the sense that one cannot obtain polynomials of higher degree for a fixed number of matrix multiplications.
Since matrix multiplications are asymptotically much more expensive than matrix additions, we can consider a larger class of polynomials obtained by taking linear combinations of previously computed matrices.
This will increase the number of degrees of freedom at a negligible cost.
We now introduce the \emph{degree-optimal polynomials}, which are constructed recursively by forming, at each step, the product of two linear combinations of all previously computed matrices.
These reach the optimal degree, and parameterize all polynomials that can be evaluated  with at most $m$ non-scalar multiplications.

Let us define the recursion
\begin{subequations}\label{eq:BBC_A}
  \begin{align}
    B_1&=I,\label{eq:BBC_A_a}\\
    B_2&=A,\\
    B_3&=(x_{a,1,1}B_1+x_{a,1,2}B_2)(x_{b,1,1}B_1+x_{b,1,2}B_2),   \\
    B_4&=(x_{a,2,1}B_1+x_{a,2,2}B_2+x_{a,2,3}B_3)(x_{b,2,1}B_1+x_{b,2,2}B_2+x_{b,2,3}B_3),\\
    B_5&=(x_{a,3,1}B_1+x_{a,3,2}B_2+x_{a,3,3}B_3+x_{a,3,4}B_4)(x_{b,3,1}B_1+x_{b,3,2}B_2+x_{b,3,3}B_3+x_{a,3,4}B_4),\\
       &\vdots  \notag\\
    B_{m+2}&=(x_{a,m,1}B_1+\cdots+x_{a,m,m+1}B_{m+1})(x_{b,m,1}B_1+\cdots+ x_{b,m,m+1}B_{m+1}).\label{eq:BBC_A_e}
  \end{align}
  As a polynomial in $A$, the matrix $B_i$ has degree 0 for $i=1$ and degree $2^{i-2}$ for $i=2$, \ldots, $m+2$, thus the linear combination
  \begin{equation}\label{eq:BBC_p}
    p(A)=y_1B_1+y_2B_2+\cdots+y_{m+2}B_{m+2},
  \end{equation}
  is a polynomial of degree at most $2^m$.
\end{subequations}
Determining the coefficients in~\eqref{eq:BBC_A} is the most challenging aspect
of this approach: in a sense, there is at once too many and too few of them.
\begin{itemize}
\item The coefficients $x_{a,1,1},\ldots,x_{b,m,m+1}$ and $y_1,\ldots,y_{m+1}$ parameterize only a subset of the vector space of polynomials of degree at most~$2^m$.
  It is easy to see that~\eqref{eq:BBC_A} has only $m^2 + 4m + 2$ parameters, which is strictly less than the $2^m$ coefficients of the polynomial for $m \ge 6$.
\item There is redundancy in the coefficients: setting $x_{a,1,1}=1$, for example, does not change the set of parameterized polynomials, and neither does setting $B_3=A^2$, i.e., $x_{a,1,1}=x_{b,1,1}=0$ and $x_{a,1,2}=x_{b,1,2}=1$, as discussed towards the end of Remark~\ref{rem:rel_not}.
\end{itemize}
The coefficients in \eqref{eq:BBC_A} can be represented in a compact form by means of two matrices of dimension $m\times (m+1)$ and one vector of length $m+2$.
We define $H_a$, $H_b$, and $y$ as
\begin{subequations}\label{eq:degopt_matrix}
  \begin{equation}\setlength{\arraycolsep}{2.5pt}
\left[\begin{array}{c|c}
    H_a & H_b
\end{array}\right]=
\left[\begin{array}{ccccc|ccccc}
  x_{a,1,1} & x_{a,1,2} &           &         &      &  x_{b,1,1} & x_{b,1,2} &           &         &      \\
  x_{a,2,1} & x_{a,2,2} & x_{a,2,3} &          &       & x_{b,2,1} & x_{b,2,2} & x_{b,2,3} &                 \\
  \vdots &             & \ddots & \ddots &       &    \vdots &             & \ddots & \ddots &       \\
  x_{a,m,1} & \cdots &  \cdots & x_{a,m,m} & x_{a,m,m+1} & x_{b,m,1} & \cdots &  \cdots & x_{b,m,m} & x_{b,m,m+1}
  \end{array}\right],
\end{equation}
\begin{equation}
y=\begin{bmatrix}y_1&\cdots&y_{m+2}\end{bmatrix}.
\end{equation}
\end{subequations}
Figure~\ref{fig:degopt_conversions} shows how several schemes for polynomial evaluation can be expressed in the form~\eqref{eq:degopt_matrix}.
\begin{figure}[t]
  \centering
  \setlength{\arraycolsep}{3.3pt}
  \subfigure[Monomial evaluation of $c_0x^0+\cdots+c_{6}x^{6}$.]{%
    \scalebox{0.81}{
      \begin{minipage}{0.60\textwidth}
        \begin{equation*}
          \left[
            \begin{array}{cccccc|cccccc}
              0& 1 &   &   &   &   &   0   & 1 &   &   & \\
              0& 0 & 1 &   &   &   &   0   & 1 & 0 &   & \\
              0& 0 & 0 & 1 &   &   &   0   & 1 & 0 & 0 &  \\
              0& 0 & 0 & 0 & 1 &   &   0   & 1 & 0 & 0 &0 \\
              0& 0 & 0 & 0 & 0 & 1 &   0   & 1 & 0 & 0 &0 &0\\
            \end{array}
          \right]
        \end{equation*}\vspace{3pt}
        \begin{equation*}
          y=\begin{bmatrix}c_0&c_1&c_2&c_3&c_4&c_5&c_6\end{bmatrix}
        \end{equation*}\vspace{3pt}
      \end{minipage}
    }\label{fig:monomial_degopt}}
  \subfigure[Horner evaluation of $c_0x^0+\cdots+c_{6}x^{6}$.]{%
    \scalebox{0.81}{
      \begin{minipage}{0.60\textwidth}
        \begin{equation*}
          \left[
            \begin{array}{cccccc|cccccc}
              c_5&c_6 &   &   &     & &   0   & 1 &   &   & \\
              c_4& 0 & 1 &   &     & &   0   & 1 & 0 &   & \\
              c_3& 0 & 0 & 1 &     & &   0   & 1 & 0 & 0 &  \\
              c_2& 0 & 0 & 0 & 1   & &   0   & 1 & 0 & 0 &0 \\
              c_1& 0 & 0 & 0 & 0   &1&   0   & 1 & 0 & 0 &0 & 0 \\
            \end{array}
          \right]
        \end{equation*}\vspace{3pt}
        \begin{equation*}
          y=\begin{bmatrix}c_0& 0 & 0 & 0 & 0 & 0 & 1  \end{bmatrix}
        \end{equation*}\vspace{3pt}
      \end{minipage}}}
  \subfigure[Paterson--Stockmeyer
  evaluation of $c_0x^0\!+\!\cdots\!+\!c_{11}x^{11}$.]{%
    \scalebox{0.81}{
      \begin{minipage}{0.60\textwidth}
        \begin{equation*}\phantom{\left[
              \begin{array}{c} 0 \\ 0 \\ 0 \\ 0 \\ 0 \\ 0 \\ \end{array}
            \right]}
          \left[
            \begin{array}{cccccc|cccccc}
              0 & 1  &  &  &     &  &   0   & 1 &   &   & \\
              0 & 0 & 1 &   &    &  &   0   & 1 & 0 &   & \\
              0 & 0 & 0 & 1 &    &  &   0   & 1 & 0 & 0 &\\
              0 & 0 & 0 & 0 & 1  &  & c_7   &c_9&c_{10}&c_{11}&0 \\
              0 & 0 & 0 & 0 & 1  &0 & c_4   &c_5&c_6& c_7& 0& 1\\
            \end{array}
          \right]
          \phantom{
            \left[
              \begin{array}{c} 0 \\ 0 \\ 0 \\ 0 \\ 0 \\ 0 \\ \end{array}
            \right]}
        \end{equation*}\vspace{3pt}
        \begin{equation*}
          y=\begin{bmatrix}c_0&c_1&c_2&  c_3 & 0&  0&1  \end{bmatrix}
        \end{equation*}\vspace{3pt}
      \end{minipage}
    }\label{fig:ps_degopt}}
  \subfigure[Truncated exponential using {\cite[Equations~(26)-(28)]{sid19}}.]{%
    \scalebox{0.81}{
      \begin{minipage}{0.60\textwidth}
        \begin{equation*}
          \left[
            \begin{array}{ccccccc|ccccccc}
              0 & 1  &  &  &     &  &  &   0   & 1 &   &   & \\  
              0 & 0 & 1 &   &    &  &  &   0   & 1 & 0 &   & \\  
              0 & 0 & 0 & 1 &    &  &  &   0   & 1 & 0 & 0 &\\   
              0 & 0 & 0 & 0 & 1  &  &  &   0   &c_4&c_3&c_2&c_1\\   
              0 &c_8&c_7&c_6&c_5 &1 &   &   0&0&c_{11}&c_{10}&c_{9}&1 \\  
              0 &c_{16}&c_{15} &c_{14} &c_{13}&c_{12}&1& 0& c_{20}   &c_{19}&c_{18}&c_{17}& 1& 0\\
            \end{array}
          \right]
        \end{equation*}\vspace{3pt}
        \begin{equation*}
          y=\begin{bmatrix}1&1&c_{23}&  c_{22} & c_{21}&  0&0 &1  \end{bmatrix}
        \end{equation*}\vspace{3pt}
      \end{minipage}
    }\label{fig:sid_degopt}}
  \caption{Various polynomial evaluation schemes
    expressed in degree-optimal form using the notation in~\eqref{eq:degopt_matrix}.}
  \label{fig:degopt_conversions}
\end{figure}

Some evaluation schemes in the literature~\cite{bbc19, sast18}
can be expressed in a degree-optimal form where certain coefficients
are fixed to zero or one---in fact \eqref{eq:BBC_A} is mentioned, though not used  extensively, in~\cite{bbc19}.
In these references, the free coefficients are computed
from algebraic system
of equations so that the resulting polynomial matches
the coefficients of a Taylor expansion of the function of interest.
The algebraic system is solved symbolically, but the solver may fail
to find a solution if the ``wrong'' coefficients are fixed, and
different choices of fixed parameters can lead
to an infinity of solutions, one of which has to
be selected in practice~\cite{bbc19}.
The authors also point out that approaches that are more economical than those obtained by
their particular choice of free and fixed parameters may exists.
As can be seen in Table~\ref{tbl:deg},
both approaches lead to polynomials of degree
far lower than the maximum obtainable and considerably below
the number of degrees of freedom
in the scheme~\eqref{eq:BBC_A}. New approaches
to determine the coefficients using symbolic computation have recently
been proposed~\cite{saib21}.

\begin{remark}[Related notation]\label{rem:rel_not}
  Some literature focusing on minimizing the number of matrix
  multiplications~\cite{sast18,sid19,siap19} uses a slightly different
  compact form to express multiplication-efficient polynomial evaluation.
  The authors prefer the notation $y_{ks}$ where $s$ is the highest power
  of the input matrix that has to be computed and $k$ is a level
  parameter such that the number of matrix multiplications is $k+s$.
  These compact forms can always be converted to a degree-optimal form \eqref{eq:BBC_A}.
  We exemplify this with the expressions in~\cite[Example~3.1]{sast18}:
  \begin{subequations}\label{eq:sastre_example31_org}
  \begin{align}
    y_{02}(x)&=x^2(c_4x^2+c_3x),   \\
    y_{12}(x)&=(y_{02}(x)+d_2x^2+d_1x)(y_{02}(x)+e_2x^2)+
    e_0y_{02}(x)+f_2x^2+f_1x+f_0.
  \end{align}
  \end{subequations}
  For any polynomial of degree 8 with non-zero leading coefficient, there exists  coefficients
  $c_4$, $c_3$, $d_2$, $d_1$, $e_2$, $e_0$, $f_2$, $f_1$, and $f_0$ such that $y_{12}(x)$ is an evaluation scheme of the polynomial~\cite{sast18}.
  In the degree-optimal form, we have one row per multiplication,
  and a final summation of all computed matrices. The coefficients
  of the degree-optimal form  for \eqref{eq:sastre_example31_org} are
\begin{equation}\label{eq:sastre_example31}\setlength{\arraycolsep}{4.25pt}
\left[\begin{array}{cccc|cccc}
  x_{a,1,1} & x_{a,1,2} &           &         & x_{b,1,1} & x_{b,1,2} &           \\
  x_{a,2,1} & x_{a,2,2} & x_{a,2,3} &           &x_{b,2,1} & x_{b,2,2} & x_{b,2,3} &           \\
  x_{a,3,1} & x_{a,3,2} & x_{a,3,3} & x_{a,3,4} &x_{b,3,1} & x_{b,3,2} & x_{b,3,3} & x_{b,3,4} \\
\end{array}\right]=
\left[\begin{array}{cccc|cccc}
      0& 1 &   &      &   0   & 1 &   &    \\
      0& 0 & 1 &      &   0   &c_3&c_4&    \\
      0& d_1 & d_2 & 1&   0   & 0 &e_2&  1  \\
    \end{array}
    \right]
\end{equation}
with
\begin{equation}\label{eq:sastre_example31_y}
\begin{bmatrix}y_1&y_2&y_3&y_4&y_5\end{bmatrix}=
\begin{bmatrix}f_0&f_1&f_2&e_0&1\end{bmatrix}.
\end{equation}
Similarly, with the notation \eqref{eq:degopt_matrix}, we can express \cite[Equation~(34)-(35)]{sast18} for $s=3$ as
\begin{equation}\label{eq:second-example-sastre}
\left[\begin{array}{ccccc|ccccc}
    0& 1 &   &   &      &   0   & 1 &   &   & \\
    0& 0 & 1 &   &      &   0   & 1 & 0 &   & \\
    0& 0 & 0 & 1 &      &   0   &c_4&c_5&c_6&  \\
    0&d_1&d_2&d_3& 1    &   0   & 0 &e_2&e_3&1 \\
 \end{array}
    \right]
\end{equation}
and
\[
\begin{bmatrix}y_1&y_2&y_3&y_4&y_5&y_6\end{bmatrix}=
\begin{bmatrix}f_0&f_1&f_2&f_3&e_0&1  \end{bmatrix}.
\]
The referenced work contains an extra addition of previous iterates after
each multiplication.
In order to transform this scheme to degree-optimal form, these terms can always be compensated for by merging into appropriate coefficients in subsequent rows in the degree-optimal form.

The function \verb#graph_sastre_yks_degopt# automatically converts the coefficients in the scheme~\cite[Equations~(62)-(65)]{sast18} for the evaluation of $y_{ks}$ to degree-optimal form.
We wish to stress that there is a redundancy in the degree-optimal form,
and 
also in the $y_{ks}$-form proposed in \cite{sast18}.
The first row of~\eqref{eq:sastre_example31} and~\eqref{eq:second-example-sastre}, for example, can be set to
\begin{equation}\label{eq:firstrow}
\left[\begin{array}{cc|cc}x_{a,1,1} & x_{a,1,2} & x_{b,1,1} & x_{b,1,2} \end{array}\right]
=
\left[\begin{array}{cc|cc}0 & 1 & 0 & 1 \end{array}\right],
\end{equation}
without loss of generality.
Our software is based on the degree-optimal polynomial
with arbitrary coefficients,
and does not explicitly eliminate this and other redundancies
that are instead removed in the $y_{ks}$-form.
This is mostly due to implementation considerations:
$H_a$ and $H_b$ having the same structure leads to more compact code.
It is straightforward to impose constraints such as~\eqref{eq:firstrow}
when optimizing the coefficients of the graph,
since we have the freedom to choose the parameters over which to optimize.
\end{remark}

The general scheme \eqref{eq:BBC_A} can clearly be phrased as
a computational graph. We provide graph generators that use
the coefficients tabulated in previous works for certain cases.
In general, we propose to approach the choice
the coefficients in~\eqref{eq:BBC_A} as an optimization problem,
as discussed in details in Section~\ref{sec:opt}.
As we will see in the experiments, determining
the coefficients by optimization can bring several advantages.
In this setting, procedures such as those in~\cite{bbc19} and~\cite{sast18}
can be used to warm start the optimizer.

\begin{table}[t]
  \begin{center}
    \caption{Highest-degree polynomial that can be evaluated with a given number
      of multiplications using different schemes.
      The asterisk ${}^*$ denotes the degree reported in the referencesd paper and signals that the corresponding coefficients were not disclosed.
      Sastre's approach~\cite{sast18} is already optimal for $m=0$, $1$, $2$, and $3$ since the degree of the polynomial is the same as the degree of \eqref{eq:BBC_p}. The approach in \cite{sid19}
      is very close to optimal for $m=4$. The plus sign ${}^{+}$ marks the highest coefficient in the Taylor series expansion that is matched by the polynomial, the approximating polynomial is higher. Sastre \cite{sast18} refers to the variants in
      \cite[Section~3.1]{sast18} and
      \cite[Section~5]{sast18}.
      \label{tbl:deg}}
    \setlength{\tabcolsep}{5.812pt}
    \begin{tabular}{@{}lcccccccc@{}}
      \toprule
      Number of matrix products $m$ & 1 & 2 & 3 & 4& 5& 6& 7&8 \\
      \midrule
      Paterson--Stockmeyer degree&                   2 & 4 & 6 & \phantom{0}9$\phantom{{}^* (15^{+})}$ & 12$\phantom{{}^+}$ & 16$\phantom{{}^+}$ &\phantom{0}20$\phantom{{}^+}$ & \phantom{0}25$\phantom{{}^+}$\\
      Sastre \cite{sast18} degree& 2 & 4 & 8 & $12^* (15^{+})$ & $16^*$& $20^*$ & $\phantom{0}25^*$ & \phantom{0}30$\phantom{{}^+}$ \\
      BBC \cite{bbc19} degree&     2 & 4 & 8 & 12$\phantom{{}^* (15^{+})}$ & 18$\phantom{{}^+}$ & $22^*$ & \\
      SID \cite{sid19} degree&     2&  4 & 8 &$15^{+}\phantom{(15^{+})}$&$21^{+}$ &24$\phantom{{}^+}$&\phantom{0}30$\phantom{{}^+}$\\
      \midrule
     Degree of \eqref{eq:BBC_p}: $2^{m}$ & 2& 4 & 8 & 16$\phantom{{}^* (15^{+})}$ & 32$\phantom{{}^+}$ & 64$\phantom{{}^+}$ & 128$\phantom{{}^+}$&256$\phantom{{}^+}$\\
      Number of free variables in~\eqref{eq:BBC_A}: $(m+2)^2-2$ & 7& 14 & 23 & 34$\phantom{{}^* (15^{+})}$ & 47$\phantom{{}^+}$ & 62$\phantom{{}^+}$& \phantom{0}79$\phantom{{}^+}$ & \phantom{0}98$\phantom{{}^+}$\\
      \bottomrule
    \end{tabular}
  \end{center}
\end{table}

\begin{example}[degree-optimal form]\label{exmp:degopt}
  The degree-optimal polynomials are represented
  by \verb#Degopt# objects, which can be initialized with a constructor
  that takes in input the matrices in \eqref{eq:degopt_matrix}.
  A degree-optimal polynomial can be evaluated by
  using the corresponding \verb#Compgraph#, which can be generated using the
  function \verb#graph_degopt#. Here is the Paterson--Stockmeyer
  evaluation corresponding to the truncated Taylor series expansion
  of the exponential.
\begin{lstlisting}
julia> c = 1 ./ factorial.(0:11); # Truncated Taylor expansion of exp(z) at z=0
julia> HA = [0.0  1.0   0.0   0.0   0.0  0.0
             0.0  0.0   1.0   0.0   0.0  0.0
             0.0  0.0   0.0   1.0   0.0  0.0
             0.0  0.0   0.0   0.0   1.0  0.0
             0.0  0.0   0.0   0.0   1.0  0.0];
julia> HB = [0.0  1.0   0.0   0.0   0.0  0.0
             0.0  1.0   0.0   0.0   0.0  0.0
             0.0  1.0   0.0   0.0   0.0  0.0
             c[9] c[10] c[11] c[12] 0.0  0.0
             c[5] c[6]  c[7]  c[8]  1.0  0.0]
julia> y=[c[1:4]; 0; 0; 1]
julia> (g, _) = graph_degopt(Degopt(HA, HB, y));
julia> A=[1 2; 3 4] / 100;
julia> eval_graph(g, A)
2x2 Matrix{Float64}:
 1.01036    0.0205091
 0.0307637  1.04112
julia> norm(exp(A) - eval_graph(g, A))
7.042759341472802e-11
\end{lstlisting}
  Since the degree-optimal form contains many zeros in this case,
  the graph will contain many nodes corresponding to
  scalar multiplications by zero. These
  can be automatically removed by the function \verb#compress_graph!#,
  which substantially reduces the execution time
  as well as the number of allocations.
\begin{lstlisting}
julia> using BenchmarkTools;
julia> g2 = compress_graph!(deepcopy(g));
julia> @btime eval_graph($g, $A);
  1.247 ms (6271 allocations: 761.16 KiB)
julia> @btime eval_graph($g2, $A);
  99.143 μs (829 allocations: 81.49 KiB)
\end{lstlisting}
\end{example}

\section{Algorithm design by optimization}\label{sec:opt}
The graphs that can be constructed from~\ref{M1}--\ref{M3} are polynomials and
rational functions.
Let $g$ be the function underlying a graph with fixed topology. Our aim is
to determine the coefficients of the linear combinations in the graph so
that the function $g$ approximates a given function $f$ in a domain of interest.
We mostly use optimization based on the forward error in exact or high-precision arithmetic.
The analysis of the backward and round-off error are presented in
Section~\ref{sec:erroranalysis}.

Throughout this section we are going to assume that $f$ and $g$ are
analytic in a simply connected domain $\Omega\subset \CC$.
Occasionally we will use the notation $g(\,\cdot\,;c)$ to emphasize that $g$
depends on a coefficient vector $c \in \CC^K$.
The forward error is $e(A) = g(A) - f(A)$, and
for any diagonalizable matrix $A\in\CC^{n\times n}$ with spectrum
in $\Omega$ we have the bound
\begin{equation}\label{eq:diag_bound}
\norm{e(A)}
\le
\kappa(V)\max_i\abs{e(\lambda_i)}
\le
\kappa(V)\max_{z\in \Omega}\abs{e(z)},
\end{equation}
where $V \in \Cnn$ is a matrix such that $\Lambda = VAV^{-1}$ is diagonal, and
$\kappa(V) = \norm{V}\norm{V^{-1}}$ denotes the spectral condition number.
Hence, we are interested in minimizing the maximum of $|e(z)|$ over $\Omega$.
Since $f$ and $g$ are assumed to be analytic in $\Omega$, then so is $e$,
and the objective function
can be further simplified by applying the maximum modulus
principle~\cite[Section~12.1]{rudi87}, which dictates that the maximum of $|e|$
over $\Omega$ is attained on the boundary.
Therefore, if we denote the boundary of $\Omega$ with $\partial\Omega$, the task is to solve
\[
\minimize_{c\in\CC^K}\max_{z\in \partial\Omega}\abs[\big]{g(z;c)-f(z)}.
\]
The resulting optimization problem is nevertheless difficult.
We introduce a surrogate problem
by replacing the $L^\infty(\partial\Omega)$-norm in the objective function
with the $L^2(\partial\Omega)$-norm
and discretizing the boundary $\partial\Omega$ obtaining the points $z_1,z_2,\dots,z_N$.
Thus, the surrogate problem is
\begin{equation}\label{eq:NLSq}
\minimize_{c\in\CC^K}
\norm{r}^2,
\end{equation}
where $r$ is the vector of residuals
\begin{equation}\label{eq:r}
r = \begin{bmatrix}
g(z_1;c) - f(z_1)\\
g(z_2;c) - f(z_2)\\
\vdots \\
g(z_N;c) - f(z_N)
\end{bmatrix}.
\end{equation}
The optimization problem \eqref{eq:NLSq} is a nonlinear least square problem
which can be tackled
with methods such as the Gauss--Newton algorithm~\cite{bjor96};
see \cite{sbl12} for an extension to the complex case.
The method updates the current guess for $c$ with a vector $\delta$ that
is the solution to the linear least squares problem
\begin{equation}\label{eq:LLSq}
\arg\min_{\delta\in\CC^K} \| J\delta - r\|,
\end{equation}
where $J$ is the Jacobian of $g$ with respect to the coefficient vector $c$, which can be written as
\begin{equation}\label{eq:J}
J = \begin{bmatrix}
\partial g(z_1;c)/\partial c_1 & \partial g(z_1;c)/\partial c_2 & \dots & \partial g(z_1;c)/\partial c_K\\
\partial g(z_2;c)/\partial c_1 & \partial g(z_2;c)/\partial c_2 & \dots & \partial g(z_2;c)/\partial c_K\\
\vdots & \vdots & \ddots & \vdots \\
\partial g(z_N;c)/\partial c_1 & \partial g(z_N;c)/\partial c_2 & \dots & \partial g(z_N;c)/\partial c_K
\end{bmatrix}.
\end{equation}
In our setting, the Jacobian can be computed in a manner that resembles backpropagation. This technique is commonly used to compute
derivatives when training neural networks \cite{rhw86,gbc16},
which can themselves
be regarded as computational graphs \cite[Section~6.5]{gbc16}.
Backpropagation applies the chain rule to the operations in the computational
graph, and by exploiting the structure of the DAG
we can derive explicit rules for the derivatives.
The derivative of~\eqref{eq:fundop} with respect to the coefficient $c_k$ is
\newsavebox{\casescommonpartder}
\begin{subequations}\label{eq:derprop}
  \begin{align}
    \sbox{\casescommonpartder}{$\displaystyle \frac{\partial Z_i(x)}{\partial c_k} =\left\{%
    \begin{array}{@{}l@{}}%
      \vphantom{\displaystyle c_i^{(1)} \frac{\partial Z_\leftp(x)}{\partial c_k} + c_i^{(2)} \frac{\partial Z_\rightp(x)}{\partial c_k}}\\[2pt]%
      \vphantom{\displaystyle \frac{\partial  Z_\leftp(x)}{\partial c_k}Z_\rightp(x) + Z_\leftp(x)\frac{\partial  Z_\rightp(x)}{\partial c_k}}\\[2pt]%
      \vphantom{\displaystyle -Z_\leftp(x)^{-2}\frac{\partial Z_\leftp(x)}{\partial c_k} Z_\rightp(x) +Z_\leftp(x)^{-1}\frac{\partial Z_\rightp(x)}{\partial c_k}}\\[2pt]%
    \end{array}%
    \right.\kern-\nulldelimiterspace$}
    \raisebox{-.66\ht\casescommonpartder}[0pt][0pt]{\usebox{\casescommonpartder}}%
    \parbox{6.375cm}{$\displaystyle c_i^{(1)} \frac{\partial Z_\leftp(x)}{\partial c_k} + c_i^{(2)} \frac{\partial Z_\rightp(x)}{\partial c_k},$} \qquad p_i=+,\label{eq:lincomb_der}\\
    \parbox{6.375cm}{$\displaystyle \frac{\partial  Z_\leftp(x)}{\partial c_k}Z_\rightp(x) + Z_\leftp(x)\frac{\partial  Z_\rightp(x)}{\partial c_k},$}\qquad p_i=*, \label{eq:mult_der}\\
    \parbox{6.375cm}{$\displaystyle -Z_\leftp(x)^{-2}\frac{\partial Z_\leftp(x)}{\partial c_k} Z_\rightp(x) +Z_\leftp(x)^{-1}\frac{\partial Z_\rightp(x)}{\partial c_k},$}\qquad p_i=\backslash, \label{eq:ldiv_der}
  \end{align}
  where $\ell = e_i^{(1)}$ and $r = e_i^{(2)}$ as in~\eqref{eq:fundop}.
  The rules~\eqref{eq:derprop} can be applied recursively until a linear combination node $s$ such that $c_k=c_s^{(1)}$ or $c_k=c_s^{(2)}$
  is found. In this case, one has
\begin{equation}
 \frac{\partial Z_s(x)}{\partial c_k} = Z_k(x). \label{eq:base_der}
\end{equation}
\end{subequations}
The propagation formulae~\eqref{eq:derprop} can also be used to implement the forward evaluation of the partial derivatives. To evaluate the derivative of $g$ with respect to $c_k$ at a point $x$
one can start from node $s$, using \eqref{eq:base_der}, and then from that node follow the parent list of the graph according to the topological ordering. Visiting the nodes in this order guarantees that whenever one of the formulae \eqref{eq:lincomb_der}-\eqref{eq:ldiv_der} is evaluated, each partial derivative on the right-hand side of \eqref{eq:derprop} has already been computed or is zero. The computation can also be carried out simultaneously for all the points $z_1,z_2,\dots,z_N$.

\begin{example}[Jacobian evaluation]\label{exmp:Jac}
  The function \verb#eval_jac# evaluates the Jacobian of a graph using the method described above.
  We illustrate the process using as model the graph for the degree-5 Taylor approximant to the exponential in monomial form.
  Evaluating this approximant requires only 4 multiplications.
\begin{lstlisting}
julia> c = 1 ./ factorial.(0:5); # Truncated Taylor expansion of exp(z) at z=0
julia> (graph, cref) = graph_monomial(c); # Exponential coeffs
julia> discr = 0.45*exp.(1im * range(0, 2π, length=200)); # Discretization
julia> J = eval_jac(graph, discr, cref);
julia> size(J)
(200, 6)
julia> svdvals(J)
6-element Vector{Float64}:
 14.142189931772608
  6.363885389264312
  2.863711391309838
  1.2886540714299903
  0.579886987450398
  0.2609452239018225
\end{lstlisting}
The second output of the graph generator, \verb#cref#, is a vector of
tuples \verb#(#\verb#Symbol#, \verb#Int#\verb#)# that references a subset
of the coefficients in the graph.
In this case, \verb#cref# contains the 6 coefficients of the monomial.
\begin{lstlisting}
julia> [cref  get_coeffs(graph, cref)]
6×2 Matrix{Any}:
 (:P2, 1)  1.0
 (:P2, 2)  1.0
 (:P3, 2)  0.5
 (:P4, 2)  0.166667
 (:P5, 2)  0.0416667
 (:P6, 2)  0.00833333
\end{lstlisting}
If the same graph is embedded in degree-optimal form, then \verb#cref#
contains all the coefficients listed in~\eqref{eq:degopt_matrix}.
\begin{lstlisting}
julia> (graph, cref) = graph_monomial_degopt(c);
julia> J = eval_jac(graph, discr, cref);
julia> size(J)
(200, 34)
julia> svdvals(J)[1:10]
10-element Vector{Float64}:
 14.142241129421516
  7.886645240829699
  3.6166303096103856
  1.347686765969578
  0.590394286817795
  0.26144313322720497
  0.005193701527333443
  0.00046141908865473055
  0.0001981485633788004
  8.380655750040214e-16
\end{lstlisting}
The flexibility of the degree-optimal form is evident, as 4 multiplications are parameterized by 34 coefficients. The redundancies are also apparent:
in this example only 9 singular values of the Jacobian are numerically larger than zero.
These 9 degrees of freedom are more than the 6 available in the
corresponding classical polynomial evaluation schemes.
Note that the number of zero singular values is a local property, which depends on the graph as well as on the values in \verb#c#.
No global conclusions, for example about which polynomials can be represented,
can be drawn from it.
\end{example}

\begin{algorithm2e}[t]
\caption{Gauss--Newton algorithm for the graph nonlinear least squares \eqref{eq:NLSq}}\label{alg:GN}
\SetKwInput{KwIn}{input}
\SetKwInput{KwOut}{output}
\KwIn{A graph, the coefficients $c_0$ corresponding to $g(\cdot\,;c_0)$, the step length $\gamma>0$.}
\KwOut{A coefficient vector $c$ such that $g(\cdot\,;c)\approx f$, if the iterations converge.}
\BlankLine
\For{$k=1,2,3,\dots$until convergence}{
Evaluate residual vector $r$ from \eqref{eq:r}.\;
Evaluate the Jacobian $J$ from \eqref{eq:J} by using the rules in \eqref{eq:derprop}.\;
Solve \eqref{eq:LLSq} for $\delta_k$ using the SVD and the pseudoinverse with a drop tolerance $d_k$. \;
Update $c_{k+1} = c_k - \gamma \delta_k$.\;
}
\end{algorithm2e}

The general idea of the Gauss--Newton algorithm adapted to the
minimization problem~\eqref{eq:NLSq} is summarized in Algorithm~\ref{alg:GN}.
Currently, our package provides a custom implementation of the Gauss--Newton
algorithm which gives good control of the damping parameters and makes it
easy to use high-precision arithmetic, via the \verb#BigFloat# data type, and to work with
real coefficients, see Remark~\ref{rem:real-NLSq}.
Optimizers from the Julia ecosystem do not provide the
features needed for our purposes, but
may be integrated in the future.
That would allow for the use of more advanced routines for
nonlinear least squares problems, such as, for example,
the Levenberg--Marquardt algorithm or general trust region
methods~\cite{bjor96}.
\begin{remark}[Real-valued coefficients]\label{rem:real-NLSq}
The nonlinear least squares problem~\eqref{eq:NLSq} is a complex-valued problem. However, in many cases it is desirable to require that ${c\in\RR^K}$. The problem can then be turned into a real-valued nonlinear least squares problem by reformulating the objective function as
\[
\minimize_{c\in\RR^K}
\left\| \begin{bmatrix}
\re(r)\\
\im(r)
\end{bmatrix}\right\|^2,
\]
where $r$ is given by \eqref{eq:r}.
The corresponding linear least squares (sub-)problem, cf.~\eqref{eq:LLSq}, becomes
\[
\arg\min_{\delta\in\RR^K} \left\| \begin{bmatrix}
\re(J)\\\im(J)
\end{bmatrix} \delta -  \begin{bmatrix}
\re(r)\\\im(r)
\end{bmatrix}\right\|,
\]
as described in~\cite{sbl12}.
\end{remark}
\begin{remark}[Relative error]
The minimization problem~\eqref{eq:NLSq} is derived from the absolute error. However, we often want to minimize the relative error $\|e(A)\|/\|f(A)\|$. If $f$ has no roots inside the domain $\Omega$, then the derivation is analogous and the result is that the Jacobian $J$ and the residual vector $r$ in the linear least squares problem~\eqref{eq:LLSq} are diagonally scaled with $D=\diag(1/f(z_i))$, i.e., the linear least squares (sub-)problem in the Gauss--Newton algorithm is $\arg\min_\delta \| DJ\delta - Dr\|$. If $f$ has roots in the domain, however, then the maximum modulus principle is not applicable to the relative error, hence we must sometimes resort to the absolute error.
\end{remark}
\begin{example}[Design by optimization]\label{exmp:optim}
We illustrate how the optimization can be used to find the coefficients of a degree-optimal polynomial of the form~\eqref{eq:BBC_A} with 4 multiplications that approximates the exponential.
This example builds on Example~\ref{exmp:Jac}, and is a simplified version
of the more extensive cases reported in Section~\ref{sec:use-cases}.
As starting guess we use the coefficients of the degree-5 Taylor approximant embedded in degree-optimal form, which has the same structure as the approximant in Figure~\ref{fig:monomial_degopt}.
\begin{lstlisting}
julia> (graph, cref) = graph_monomial_degopt(c); # Initial guess
julia> f = exp;
julia> A = randn(100, 100) / 40;
julia> fA = f(A);
julia> norm(eval_graph(graph, A) - fA) / norm(fA) # Test initial guess
3.521697764659818e-7
\end{lstlisting}
The second output in the graph generation, \verb#cref#, is described in Example~\ref{exmp:Jac},
and represents the coefficients we use in the optimization.
We apply Gauss--Newton to find an improved graph.
\begin{lstlisting}
julia> graphb = big(graph); # Do optimization in high precision
julia> discr = 0.45 * exp.(1im*range(0,2big(π), length=200)); # Discretization
julia> opt_gauss_newton!(graphb, f, discr,  # Run optimization
                         errtype = :relerr, # Relative error
                         stoptol = 4e-15,
                         cref = cref, # Explicitly giving target variables
                         linlsqr = :real_svd, # Force real coefficients
                         droptol = 1e-15,
                         logger = 0 # Set to 1 to see iterations
                         );
julia> graph = Compgraph(Float64, graphb); # Round to Float64
julia> norm(eval_graph(graph, A) - fA) / norm(fA) # Test result
3.3790722357852773e-16
\end{lstlisting}
The relative error is of the order of the unit round-off, an improvement
of 9 orders of magnitude over the truncated Taylor approximant.
This illustrates the potential of the optimization approach, but further
analysis is needed in order to establish whether a graph found by optimization
will be useful in practice.
In Section~\ref{sec:bwerrexp} we provide tools for a detailed error analysis,
and study the accuracy of \verb#graphb# in Example~\ref{exmp:bwerr}.
A more complete design procedure is outlined in Section~\ref{sec:use-cases}, where we also
present improved graphs.
\end{example}

\begin{remark}[Jordan form and generalization of diagonal bound~\eqref{eq:diag_bound}]
If $A$ is defective, then defining $f(A)$ requires the ability to compute some of the derivatives of $f$ at some eigenvalues of~$A$, as can be seen from the definitions via Jordan canonical form~\cite[Definition~1.2]{high08} or via Hermite interpolation~\cite[Definition~1.3]{high08}.
Therefore, in order to minimize $\norm{e(A)}$ it might be desirable to ensure that the derivatives of $g$ approximate those of $f$.
Since $f$ and $g$ are analytic on $\Omega$, so is the error $e$, and we can thus bound its derivatives.

To do this we consider $\Omega$ being embedded in a slightly larger region $\Omega_\rho$.
If $D(a,\rho)$ denote the disk of radius $\rho$ centered in $a$, then
$\Omega_\rho$ is the smallest region such that for every point $a\in\Omega$
one has $D(a, \rho) \subset \Omega_\rho$.
We assume that $f$ and $g$ are analytic in $\Omega_\rho$.
For reasons that will become clear later, it is preferable to set $\rho\geq1$.
Since $D(a,\rho)\subset\Omega_\rho$ for all points $a\in\Omega$, by Cauchy's estimates \cite[Theorem~10.26]{rudi87} we have that
\[
\abs[\big]{g^{(k)}(a)-f^{(k)}(a)}\leq \frac{k! M_{a,\rho}}{\rho^k},
\]
where $M_{a,\rho}$ is a constant such that $\abs[\big]{g(z)-f(z)}\leq M_{a,\rho}$ for all $z\in D(a,\rho)$.
Naturally for any $a\in\Omega$ we have $M_{a,\rho}\leq M_\rho$ where $M_{\rho}$ is a constant such that $|g(z)-f(z)|\leq M_{\rho}$ for all $z\in \Omega_\rho$. Hence, the derivative bounds can be applied uniformly giving
\[
\abs[\big]{g^{(k)}(z)-f^{(k)}(z)}\leq \frac{k! M_{\rho}}{\rho^k},
\]
for all $z\in\Omega$.
Thus, if $g$ approximates $f$ well in the slightly larger domain $\Omega_\rho$,
then we can bound how well the
derivatives of $g$ approximate those of $f$ on $\Omega$.
In the context of matrix functions, this implies that for a $k \times k$ Jordan block $J_k$ one has that
\[
\norm[\big]{g(J_k)-f(J_k)}
\leq
\norm[\big]{\,\abs[\big]{g(J_k)-f(J_k)}\,}
\leq
M_\rho \norm{P_k(\rho)},
\]
where $|\cdot|$ denotes the entry-wise absolute value and  $P_k(\rho)$ is the matrix
\[
P_k(\rho) = \begin{bmatrix}
1 & \frac{1}{\rho} & \frac{1}{\rho^2} & \dots & \frac{1}{\rho^{k-1}}\\
  & 1 & \frac{1}{\rho} & \ddots & \vdots \\
  &   & 1 & \ddots & \frac{1}{\rho^2} \\
  &   &   & \ddots & \frac{1}{\rho}\\
  &   &   &  & 1
\end{bmatrix}.
\]
We can conclude that
\begin{equation*}
  \|e(A)\|\leq \kappa(S) \|P_{k'}(\rho)\| \max_{z\in \Omega_\rho}|e(z)|,
\end{equation*}
where $k'$ is the size of the largest Jordan block of $A$, and $S$ is a transformation matrix such that $J = SAS^{-1}$ is a Jordan canonical form of $A$.
\end{remark}

\begin{remark}[Crouzeix's conjecture and generalization of diagonal bound~\eqref{eq:diag_bound}]
  Let $\Omega$ be a smooth, bounded, and
  convex domain that contains the field of values
  \begin{equation*}
    W(A) =\{x^*Ax: x\in\CC^n, \|x\|=1\}.
  \end{equation*}
  By~\cite[Theorem~3.1]{crpa17}, if $e$ is analytic in $\Omega$ and continuous on its boundary, then
\[
  \norm{e(A)}\leq \Bigl(1+\sqrt{2}\Bigr)\sup_{z\in \Omega} \abs{e(z)}.
\]
Let $\Omega=D(0,r)$ be the disk with radius $r$ centered at the origin, and let us assume that $f$ and $g$ are analytic in $D(0,r)$,
which implies that so is the error $e$.
Consider the numerical radius $w(A)=\sup\{|z|:z\in W(A)\}$.
For any matrix $A$ such that $\|A\|\leq r$ the numerical radius satisfies $w(A)\leq r$
\cite[Problem~1.5.24(g)]{hj:topics}, and thus, the numerical range satisfies $W(A)\subset D(0,r)$.
Hence, for all matrices $\|A\|\leq r$ we have that
\[
\|e(A)\|\leq \Bigl(1+\sqrt{2}\Bigr)\sup_{z\in D(0,r)} |e(z)|.
\]
\end{remark}

\section{Efficiency and reproducibility features}
\label{sec:eff-rep}
The package is written in Julia but provides a wide range of facilities
to ensure that users of other languages can benefit from our findings.
Therefore, we have implemented functions for saving and loading
graphs in a language-independent format and for
generating multi-language code directly from a graph object.

\subsection{Code generation}\label{sec:codegen}
The package allows the user to generate fast and memory-efficient code.
The code generation feature is also available for Julia,
and is in general preferred over the \verb#eval_graph# function
for efficiency.
The package supports code generation for
\begin{itemize}
\item Julia, with in-place addition to reduce the memory footprint,
\item MATLAB,
\item C with support for OpenBLAS, and
\item C with support for the Intel Math Kernel Library (MKL).
\end{itemize}

On current HPC computer architectures, memory access and allocations
have a large impact on performance. In our work, we strove
to minimize the memory footprint of the generated code by choosing a topological
ordering that (heuristically) attempts to minimize the maximum number of matrices
that the algorithms needs to store at any given time. Finding
the topological ordering that is optimal in the sense of the memory footprint
is essentially equivalent to determining the width of the
graph, which is a difficult problem in general.
For generality, the choice of topological ordering can therefore be
influenced by parameters passed to \verb#gen_code#.

\begin{example}[MATLAB code generation]\label{exmp:codegen}
  The following code generates the computational graph
  corresponding to the Paterson--Stockmeyer evaluation of
  the Taylor approximant to the cosine.
  We exploit the fact that the truncated Taylor series expansion of $\cos A$ is a polynomial in $A^2$.
\begin{lstlisting}
julia> # c is the non-zero monomial coefficients for the cosine function
julia> c = (kron(ones(5) ,[1, -1]) ./ factorial.(0:2:19));
julia> (g,_) = graph_ps(c); # Paterson--Stockmeyer evaluation
julia> rename_node!(g, :A, :A2tmp); # ``change'' input to A^2
julia> add_mult!(g, :A2tmp, :A, :A);
julia> gen_code("mycosm.m", g, funname = "mycosm", lang = LangMatlab())
\end{lstlisting}
  This generates the file:
  \begingroup
  \lstset{
    language=Matlab
  }
\begin{lstlisting}
function output = mycosm(A)
    n = size(A,1);
    I = eye(n,n);
    % Computation order:
    % A2tmp B_0_1 B_1_1 B_2_1 A2 B_0_2 B_2_2 B_1_2 A3 P2 C1 P1 C0 P0
    % Computing A2tmp with operation: mult
    A2tmp = A * A;
    % Computing B_0_1 with operation: lincomb
    coeff1 = 1.0;
    coeff2 = -0.5;
    B_0_1 = coeff1*I + coeff2*A;
    % Computing B_1_1 with operation: lincomb
    coeff1 = -0.001388888888888889;
    coeff2 = 2.48015873015873e-5;
    B_1_1 = coeff1*I + coeff2*A;
    ...
    P1 = coeff1*C1 + coeff2*B_1_2;
    % Computing C0 with operation: mult
    C0 = P1 * A3;
    % Computing P0 with operation: lincomb
    coeff1 = 1.0;
    coeff2 = 1.0;
    P0 = coeff1*C0 + coeff2*B_0_2;
    output = P0;
end
\end{lstlisting}
\endgroup
The function \verb#mycosm#
can be called directly from MATLAB:
\begin{lstlisting}
>> cosm = @(S)(expm(1i*S) + expm(-1i*S)) / 2 % Euler's formula
>> A = randn(100, 100) / 100;
>> norm(cosm(A) - mycosm(A))
ans =
   2.2390e-16
\end{lstlisting}
\end{example}
\begin{remark}[Dot fusion and in-place addition]
  In terms of performance, Julia has a competitive advantage
  for the matrix operations considered in the package.
  The reason is twofold: on the one hand the extended usage of dots
  and dot fusion~\cite{john17} reduces the number of passes through
  the data required to carry out a computation, while on the other
  the availability of in-place operations decreases the memory footprint.
  Our code generation can take full advantage of these
  two features, thus gaining considerable performance improvements.
  In fact, the dot fusion is the reason why the generated Julia code outperforms
  even C implementations in some of the experiments in
  Section~\ref{sec:use-cases}. In particular, when relying on the BLAS one can
  only compute the sum of two matrices at a time, whereas by using the dot
  fusion, Julia can add an arbitrary number of matrices with only a single pass
  through the data.
  The cosine in the previous example is used as an illustration.
  \begin{lstlisting}
julia> using BenchmarkTools
julia> # Dot fusion is enabled by default
julia> gen_code("mycosm.jl", g, funname = "mycosm", lang = LangJulia())
julia> # It can be disabled with a constructor argument
julia> use_dot_fusion = false;
julia> gen_code("mycosm_no_dot.jl",g,
                funname = "mycosm_no_dot",
                lang = LangJulia(true, true, use_dot_fusion))
julia> include("mycosm.jl");
julia> include("mycosm_no_dot.jl");
julia> A0 = randn(500, 500) / 500;
julia> A = copy(A0); @btime mycosm($A);
  20.576 ms (9 allocations: 7.63 MiB)
julia> A = copy(A0); @btime mycosm_no_dot($A);
  27.924 ms (11 allocations: 9.54 MiB)
\end{lstlisting}
\end{remark}

\subsection{Storing a graph in a file}\label{sect:cgr}

\begin{figure}[t]
  \centering
\begin{Verbatim}[frame=single]
% # Representation of a computation graph
% # Created: 2021-05-25T11:48:01.720 by user jarl

graph_coeff_type="Float64";

A2tmp=A*A;
coeff1=1.0;
coeff2=-0.5;
B_0_1=coeff1*I+coeff2*A2tmp;
coeff1=-0.001388888888888889;
coeff2=2.48015873015873e-5;
B_1_1=coeff1*I+coeff2*A2tmp;
coeff1=2.08767569878681e-9;
coeff2=-1.1470745597729725e-11;
B_2_1=coeff1*I+coeff2*A2tmp;
...
\end{Verbatim}
  \caption{CGR file for the graph in Example~\ref{exmp:codegen}.}
  \label{fig:cgr-file}
\end{figure}

Language-independent reproducibility of our results is ensured by the
availability of functions for saving and loading \verb#Compgraph# objects.
The graphs are stored in a human-readable way that can easily be interpreted
without using Julia.

The functions \verb#export_compgraph#
and \verb#import_compgraph# read and write the graph,
including optional or generated graph metadata.
The output is saved in computational graph (CGR) format, a storage format
specifically designed to store computational graphs for matrix functions.

\begin{example}[Saving and loading]\label{exmp:cgr}
  We can export the graph generated in Section~\ref{sec:codegen}
  to a file with the command:
\begin{lstlisting}
julia> export_compgraph(g, "cosm_graph.cgr");
\end{lstlisting}
which produces the file in Figure~\ref{fig:cgr-file}.
We note that a CGR file is essentially an executable script in both
MATLAB and Julia, but we stress that this is not the recommended way
to generate code in either language.
\end{example}

\section{Error analysis}\label{sec:erroranalysis}

\subsection{Bounds on forward and backward error}\label{sec:bwerrexp}

\newcommand{\thfwd}{\ensuremath{\theta^F}}
\newcommand{\thbwd}{\ensuremath{\theta^B}}

Let us consider the problem of evaluating the function $f$ at a matrix $A$
when several graphs $g_0$, \dots, $g_\ell$ approximating $f$ are
available.
Without loss of generality, we can assume that these functions are ordered by
increasing computational cost, so that evaluating $g_{i+1}(A)$ is
computationally more expensive than evaluating $g_{i}(A)$.
In practice, we will typically also have that $\norm[\big]{g_{i+1}(A) - f(A)} <
\norm[\big]{g_i(A) - f(A)}$, but the analysis in this section does not
require this.

In order to approximate $f$ efficiently, ideally we would like to have an
inexpensive strategy to choose at runtime the smallest $i$ such that $g_i(A)$ is
guaranteed to be an acceptable approximation of~$f(A)$.
One can gauge the quality of the approximation in several ways.

In Section~\ref{sec:opt}, we described how the forward error is an effective
tool to optimize the coefficients of a computational graph for a fixed graph
topology.
In principle, the relative forward error could be used as a means to choose a
suitable approximant, as we now explain.
Our analysis is based on a result of Al-Mohy and
Higham~\cite[Theorem~4.2(a)]{alhi09a}, which we recall here as it will be
crucial to our discussion.
If the function $h : \Omega \subset \CC \to \CC$ has the series expansion
\begin{equation}
  \label{eq:hz-series}
  h(z) = \sum_{j=j_0}^\infty \beta_j z^j,
\end{equation}
then for any matrix $A$ such that
\begin{equation}
  \label{eq:alpha-def}
  \alpha_p(A) = \max \Bigl\{\norm[\big]{A^p}^{1/p},
  \norm[\big]{A^{p+1}}^{1/(p+1)}\Bigr\},
  \qquad
  p(p-1) \le j_0,
\end{equation}
is within the disk of convergence of the series, then one has that
\begin{equation}
  \label{eq:series-bound}
  \norm[\big]{h(A)} \le \sum_{j=j_0}^\infty \abs{\beta_j} \alpha_p(A)^j.
\end{equation}
We note that if the coefficients $\beta_j$ in~\eqref{eq:hz-series} are
one-signed, then the bound~\eqref{eq:series-bound} reduces to
\begin{equation*}
  \norm[\big]{h(A)} \le \abs[\big]{h\bigl(\alpha_p(A)\bigr)},
\end{equation*}
and also that~\eqref{eq:alpha-def} reduces to
\begin{equation}
  \label{eq:alpha-def-0}
   \alpha_p(A) = \max\{1, \norm{A}\}
\end{equation}
in the general case $j_0 = 0$.
The function~$\alpha_p(A)$ in~\eqref{eq:series-bound} could in principle be
replaced with alternatives that are typically more expensive to compute but
yield a sharper bound on the norm of $h(A)$. As for the series of interest we
typically have that $j_0 = 0$, we prefer the simpler result
in~\cite[Theorem~4.2(a)]{alhi09a} to, for example, \cite[Theorem~1]{sidr11} or
\cite[Theorem~1]{sird13a}, which require that $j_0$ be at least 1, or to
\cite[Theorem~3.1]{nahi18}, which reduces to~\eqref{eq:alpha-def-0} for
$j_0 = 0$.

One can then define, using the series expansion of the forward error, the series
\begin{equation*}
  E_i(z) = \sum_{j=j_0}^\infty \abs{\gamma_j} z^j,
  \quad
  \text{where}
  \quad
  e_i(z) = g_i(z) - f(z) = \sum_{j=j_0}^\infty \gamma_j z^j,
  \quad
  j_0 \in \NN,
\end{equation*}
and then pre-compute the values
\begin{equation}
  \label{eq:thfwdi}
  \thfwd_i =
  \max\bigl\{\theta \in \RR^+ : E_i\bigl(\theta\bigr) \le u \bigr\},
\end{equation}
where $\RR^+$ and $u$ denote the positive real semi-axis and the unit round-off
of the current working precision, respectively.
It is easy to see that choosing a $g_i$ such that $\alpha_p(A) \le \thfwd_i$
guarantees that the absolute forward error $\norm{e_i(A)}$ will not be larger
than~$u$.
This measure of accuracy, however, is not ideal for two reasons.
On the one hand, the use of an absolute forward error disregards the magnitude
of the matrix function one seeks to compute, although turning to relative errors is
unfeasible in most cases, as estimates of $\norm{f(A)}$ are rarely available.
On the other hand, the dependency of the forward error on the condition number
of $f$ evaluated at $A$, makes this method not completely reliable.
For these two reasons, bounds on the forward error are seldom used in
algorithms for matrix functions, one exception being arbitrary-precision
environments, in which the pre-computation of the $\thfwd_i$
in~\eqref{eq:thfwdi} would be impractical as the value of $u$ is known only at
runtime~\cite{fahi18, fahi19}.

A better alternative, originally proposed by Moler and Van Loan~\cite{mova78}
to choose the degree of a diagonal Padé approximant to the exponential for the
scaling and squaring algorithm, is based on bounds on the relative backward
error.
The essence of this approach is to consider the matrix $E$ such that $g_i(A) =
f(A+E)$ and then pick an approximant $g_i$ for which ${\norm{E} \le
u \norm{A}}$.
In order to do that, we can rely on the series expansion
\begin{equation}
  \label{eq:def-Fi}
  F_i(z) = \sum_{j=j_0}^\infty \abs{\delta_j} z^{j-1},
  \quad
  \text{where}
  \quad
  \varphi_i(z) = \sum_{j=j_0}^\infty \delta_j z^j,
\end{equation}
and the absolute backward error $\varphi_i(z)$ is defined by $g_i(z) = f(z +
\varphi_i(z))$.
If, when designing the algorithm, we compute the values
\begin{equation}
  \label{eq:thbwdi}
  \thbwd_i =
  \max\bigl\{\theta \in \RR^+ : F_i\bigl(\theta\bigr) \le u \bigr\},
\end{equation}
then the approximant $g_i$ evaluated at a matrix $A \in \Cnn$ such that
 $\alpha_p(A) \le \thbwd_i$ will guarantee that
\begin{equation*}
  \frac{\norm{\varphi_i(A)}}{\norm{A}} \le u,
\end{equation*}
which is entirely satisfactory.

In our package we provide the function
\texttt{compute\_bwd\_theta} to compute
the values $\thbwd_i$ in~\eqref{eq:thbwdi} for the exponential.
The function exploits the identity
\begin{equation}
  \label{eq:bwd-exp-series}
\varphi_i(z) = \log\bigl((e^{-z} g_i(z)-1)+1\bigr),
\end{equation}
and uses as $F_i$ the truncated series expansion of the expression on the
right-hand side of \eqref{eq:bwd-exp-series}.
Truncated series expansion are manipulated efficiently as polynomials by relying
on the functions available in the \texttt{Polynomials.jl} package.
The computation of $\thbwd_i$ is then performed using the \texttt{fzero}
function in arbitrary precision.

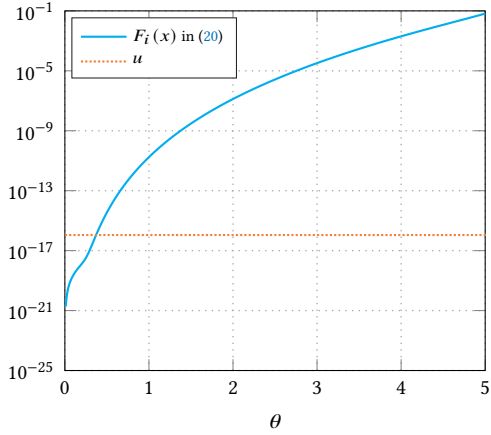
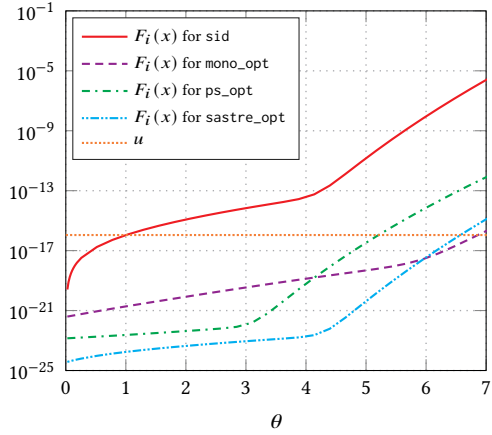
\begin{figure}[t]\centering
  \subfigure[Approximant in Example~\ref{exmp:bwerr} ($m=4$).%
  \label{fig:thbwdi_a}]{%
    \pgfplotsset{
  every axis/.append style={title style={yshift=1ex}},
  grid=major,
  grid style={line width=.5pt, draw=Black!40},
  major grid style={dotted}}

\begin{tikzpicture}[trim axis left, trim axis right]

  \pgfplotstabletranspose\loadedtable{gfx/example-bwd-err.dat}

  \begin{axis}[
    xmin=0, xmax=5,
    xtick distance={1},
    ymin=1e-25, ymax=1e-1,
    max space between ticks=25pt,
    xlabel={$\theta$},
    legend style={at={(axis cs: 0.07143,0.04)},
      anchor=north west,
      draw=white!15!black,
      legend cell align={left},
      font=\tiny,
      cells={line width=0.75pt}},
    ymode=log,
    width=2.8in,
    height=2.5in,
    cycle list name = bwd-list,
    every axis plot/.append style={thick}
    ]

    \addplot table [x index=1,y index=2] \loadedtable;
    \addplot coordinates {
      (0, 1.1e-16)
      (7, 1.1e-16)};

    \addlegendentry{$F_i(x)$ in \eqref{eq:def-Fi}}
    \addlegendentry{$u$}

  \end{axis}
\end{tikzpicture}}%
  \qquad\qquad%
  \subfigure[Approximants in Table~\ref{tbl:theta} ($m=7$, $r=6.4$).%
  \label{fig:thbwdi_b}]{\pgfplotsset{
  every axis/.append style={title style={yshift=1ex}},
  grid=major,
  grid style={line width=.5pt, draw=Black!40},
  major grid style={dotted}}

\begin{tikzpicture}[trim axis left, trim axis right]

  \begin{axis}[
    xmin=0, xmax=7,
    xtick distance={1},
    ymin=1e-25, ymax=1e-1,
    max space between ticks=25pt,
    xlabel={$\theta$},
    legend style={at={(axis cs: 0.1,0.04)},
      font=\tiny,
      anchor=north west,
      draw=white!15!black,
      legend cell align={left},
      cells={line width=0.75pt}},
    ymode=log,
    width=2.8in,
    height=2.5in,
    cycle list name = bwd-list-2,
    every axis plot/.append style={thick}
    ]
    \addplot table{gfx/bwderr_exp_sid_m7.csv}; 
    \addplot table{gfx/bwderr_exp_mono_m7_opt_rho6_4.csv}; 
    \addplot table{gfx/bwderr_exp_ps_m7_opt_rho6_4.csv}; 
    \addplot table{gfx/bwderr_exp_sastre_m7_opt_rho6_4.csv}; 
    \addplot coordinates {(0, 1.1e-16) (10, 1.1e-16)};

    \addlegendentry{$F_i(x)$ for \texttt{sid}}
    \addlegendentry{$F_i(x)$ for \texttt{mono\_opt}}
    \addlegendentry{$F_i(x)$ for \texttt{ps\_opt}}
    \addlegendentry{$F_i(x)$ for \texttt{sastre\_opt}}

    \addlegendentry{$u$}
  \end{axis}
\end{tikzpicture}}
  \caption{Series expansion~\eqref{eq:def-Fi} for graphs with and without
    optimized coefficients. The dotted orange line represents the unit round-off of
    binary64 arithmetic $u = 2^{-53} \approx 1.11 \times 10^{-16}$.
    To warm start the optimization we used for the figure on the left the graph in
    Example~\ref{exmp:optim}, which requires $m=4$ multiplications, and for
    the figure on the right several methods from the literature that
    require $m=7$ multiplications; see Section~\ref{sec:use-cases} for details.
  }
  \label{fig:bwd-err-exp}
\end{figure}

\begin{example}[Backward error function]\label{exmp:bwerr}
  We illustrate the usage of \texttt{compute\_bwd\_theta} by continuing Example~\ref{exmp:optim}.
  Because of how Julia handles types and how the \texttt{Compgraph} data structure is implemented, we can perform a symbolic evaluation of a graph to produce a \texttt{Polynomial}.
\begin{lstlisting}
julia> using Polynomials;
julia> x = Polynomial("x");
julia> p = eval_graph(graphb, x)
Polynomial(1.00000000000000000134967838007550967408645631852098 +
 1.00000000000000000432055652330969008190946510481787*x +
 0.50000000000000001033494602832231711963058471255608*x^2 +
 0.16666666666666669000159470466182272865390213333447*x^3 +
 0.04166666666666671867523152654919958881035762972403*x^4 +
 0.00833333333333344895439682913026016254992768036399*x^5 +
 0.00138888888888914582962100608532582512792424020151*x^6 +
 0.00019841269841326921065494640481934558709823865755*x^7 +
 0.00002480158730284895148988505562762560984853494876*x^8 +
 0.00000275573192499137884911178001938591470365322543*x^9 +
 0.00000027557319274655290493587144103267057249682196*x^10 +
 0.00000002505201184943177952245945443077211293385948*x^11 +
 0.00000000208628333194763098184540305832611877450755*x^12 +
 0.00000000015127436326998476974610040316369471425214*x^13 +
 0.00000000000845801512354938892869969428079487596106*x^14 +
 0.00000000000030324418961689639808749612402211552411*x^15 +
 0.00000000000000515212112436523524706264550043349319*x^16)
\end{lstlisting}
  These coefficients are close to those of the Taylor series expansion, but not
  identical, neither in arbitrary \emph{(}\verb#BigFloat#\emph{}) nor in binary64
  \emph{(}\verb#Float64#\emph{)} arithmetic.
\begin{lstlisting}
julia> a = p.coeffs;
julia> a_taylor = 1 ./ factorial.(BigFloat.(0:16));
julia> norm(abs.(a-a_taylor))
9.900938779409988e-12
julia> a_float64 = convert(Array{Float64}, a);
julia> a_taylor_float64 = convert(Array{Float64}, a_taylor);
julia> norm(abs.(a_float64 - a_taylor_float64))
9.900938779407373e-12
\end{lstlisting}

  We can approximate the bound $\thbwd_i$ in~\eqref{eq:thbwdi} by using the function \verb#compute_bwd_theta#.
  For illustration we now show how this can be computed directly by using the first 100 coefficients of the series expansion in~\eqref{eq:def-Fi} and \texttt{BigFloat} arithmetic with 1000 significant binary digits.
\begin{lstlisting}
julia> setprecision(1000); # Set number of bits in fraction of BigFloat.
julia> nterms = 100; # Order of highest term in the truncated series expansions.
julia> p = convert(Polynomial{Float64, :x}, p); # Use Float64 polynomial.
julia> # Evaluate truncated series expansion of right-hand side of (21).
julia> coeffs_exp = (-1).^(0:1:nterms) ./
                    factorial.(collect(big(0.):big(1.):nterms)); # exp(-z)
julia> expminusz = Polynomial(coeffs_exp);
julia> R = expminusz * p - 1; # Rᵢ (z) = exp(-z) p(z) - 1
julia> coeffs_log = [0; (-1) .^ (0:1:nterms-1) ./ (1:1:nterms)]; # log(x+1)
julia> logzplusone = Polynomial(coeffs_log);
julia> F = logzplusone(R); # log(Rᵢ+1).
julia> bnd_bwd_err = Polynomial(abs.(F.coeffs)) # Absolute value of coefficients
julia> # Find point where bound on relative backward error equals tolerance.
julia> using Roots
julia> bnd_rel_bwd_err(z) = abs.(bnd_bwd_err(z)) ./ abs.(z)
julia> theta_bwd = fzero(z->bnd_rel_bwd_err(z) - eps(Float64)/2, big"0.2")
0.37187598755150271424788203619606964195065108599084
\end{lstlisting}
  The function \verb#bnd_rel_bwd_err# represents a bound on the relative backward error of the function underlying \texttt{graphb} seen as an approximant to the exponential.
  The plot of the function \verb#bnd_rel_bwd_err#
  is given in Figure~\ref{fig:thbwdi_a} where we also specified
  the corresponding $\thbwd_i$. For this graph, the backward
  error will be smaller than machine precision
  for any complex scalar of modulus at most $\thbwd_i\approx 0.37$.
\end{example}

\subsection{Round-off error}

\label{sec:roundoff}

The optimization process in Section~\ref{sec:opt} and the
backward error analysis are based on exact arithmetic and rely on
high-precision arithmetic for practical purposes.
The algorithms, however, are designed to be used in floating-point arithmetic,
where rounding errors may become the most prominent cause for loss of accuracy.
The generality of the data structures in our package
allows us to obtain more realistic
estimates and bounds on the round-off errors that may occur during the
evaluation of $g(z)$.
The computational approach discussed here is based on an a posteriori error analysis of the scalar evaluation, thus it will not capture issues that pertain only to matrix evaluation, such as for example the loss of commutativity caused by round-off errors.

In particular, we will rely on a running error analysis, in which the rounding
errors are bounded automatically by the algorithm as $g(z)$ is evaluated.
Our goal will be to assign to each node $i \in V$ a value
$\delta_i$ that represents the round-off error occurred in the evaluation of
the subgraph rooted at node $i$.
In the remainder of this section a hat denotes computed quantities,
thus we will employ the notation
$\wh Z_i = Z_i(1+  \delta_i)$, where $\delta_i$ represents the relative
round-off error~\cite[Theorem~2.2]{high02}.
For conciseness, we will also denote the indices of the two parents of node
$Z_i$ by $\leftp = e_i^{(1)}$ and $\rightp = e_i^{(2)}$, as done
in~\eqref{eq:fundop}.

Note that the only two ingredients the running error analysis on the graph
requires are the value of $\delta_i$ for the leaf nodes---the input
nodes of the graph---and a rule to compute $\delta_i$ from $Z_{\leftp}$,
$Z_{\rightp}$, $\delta_{\leftp}$, and $\delta_{\rightp}$.
If $i$ is an input node, then we can set $\delta_i = 0$, which is equivalent to
saying that the input values $1$ and $z$ are represented exactly.
If $i$ is not an input node, we need to consider three cases.
\begin{itemize}
\item If $p_i = *$, then using the standard model of floating-point
  arithmetic~\cite[Equation~(2.5)]{high02} we have that
  \begin{equation*}
    \wh Z_i(1+\eta) = \wh Z_\leftp \wh Z_\rightp,
    \qquad
\abs{\eta} \le u,
\end{equation*}
and by using the fact that $Z_i=Z_\leftp Z_\rightp$,
we can solve for $\wh Z_i$ and obtain
  \begin{equation}
    \label{eq:run-err-mult}
    \wh Z_i=Z_i\frac{(1+\delta_\leftp)(1+\delta_\rightp)}{1+\eta}
  \end{equation}
  and therefore
  \[
  \delta_i=\frac{(1+\delta_\leftp)(1+\delta_\rightp)}{1+\eta}-1=\frac{\delta_\leftp+\delta_\rightp+\delta_\rightp\delta_\rightp-\eta}{1+\eta}\approx \delta_\leftp+\delta_\rightp-\eta.
  \]
\item Similarly, if $p_i = \backslash$, then we have that  $Z_i=Z_\leftp^{-1}Z_\rightp$ and
  \[
  \wh Z_i=Z_i\frac{1+\delta_\rightp}{(1+\delta_\leftp)(1+\eta)},
  \qquad
    \abs{\eta} \le u,
\]
which implies that
\begin{equation}\label{eq:run-err-ldiv}
\delta_i=\frac{\delta_\rightp-\delta_\leftp-\eta-\delta_\leftp\eta}{(1+\delta_\leftp)(1+\delta_\rightp)}\approx  \delta_\rightp-\delta_\leftp-\eta.
\end{equation}
\item Finally, if $p_i = +$, we have that
  \begin{equation*}
    \wh Z_i = \frac{\alpha \wh Z_\leftp + \beta \wh Z_\rightp}
    {1 + \eta},\qquad
    \abs{\eta} \le u.
  \end{equation*}
  We solve this equation for $\wh Z_i$ and obtain
\[
\wh Z_i=Z_i\frac{\alpha Z_\leftp (1+\delta_\leftp)+\beta Z_\rightp (1+\delta_\rightp)}{(1+\eta)Z_i},
\]
  from which we find the relative error
  \begin{equation} \label{eq:run-err-lincomb}
    \delta_i=
    \frac{\alpha Z_\leftp\delta_\leftp+\beta Z_\rightp\delta_\rightp-\eta Z_i}
    {(1+\eta)Z_i}
    \approx \frac{\alpha Z_\leftp\delta_\leftp+\beta Z_\rightp\delta_\rightp}
    {Z_i}-\eta.
  \end{equation}
\end{itemize}

The formulae for the relative error
in~\eqref{eq:run-err-mult},~\eqref{eq:run-err-ldiv},
and~\eqref{eq:run-err-lincomb} are exact,
but they are also impractical:
\begin{enumerate}
\item The quantity $\eta$ is not necessarily representable in the current working
  precision, is not readily available, and is unduly expensive to estimate.
\item The expressions for linear combinations feature exact quantities,
  which are not available during a genuine run of the algorithm.
\end{enumerate}
In our setting, the second point can be handled, in practice, by considering
estimates computed in high-precision arithmetic.

One way to make these formulae practical is to turn them into running error
bounds.
In the approximations on the right-hand side of
\eqref{eq:run-err-mult}, ~\eqref{eq:run-err-ldiv},
and~\eqref{eq:run-err-lincomb}
we have linearized the expressions by neglecting terms of the form
$\delta_\rightp^i\delta_\leftp^j\eta^k$ for $i+j+k \ge 2$.
Exact quantities can be replaced by their computed counterparts.
This substitution will not have a major effect, since the magnitude of
$\delta_i$ should be negligible compared with that of $Z_i$.
Then we can easily obtain running bounds on the relative error by replacing the
``exact'' errors $\delta_i$ in the formulae above with the approximate bounds
\begin{equation}
  \label{eq:run-err-bnd}
  \varDelta_i =
  \begin{dcases}
    \vphantom{\frac{1}{2}}
    \frac{\abs{\alpha}|\wh Z_\leftp|\varDelta_\leftp + \abs{\beta}|\wh Z_\rightp|\varDelta_\rightp}
    {\abs{\wh Z_i}} + u,\qquad &p_i = +,\\
    \vphantom{ZZZ\frac{1}{2}ZZZZZ}   \Delta_\leftp+\Delta_\rightp+u,\qquad& p_i \in \{*,\backslash\},
  \end{dcases}
\end{equation}
for which $\abs{\delta_i} \lesssim \Delta_i$.
It is important to stress that $\Delta_i$ may potentially not be a true bound
because of the linearization step.

In practice, the formula \eqref{eq:run-err-bnd} can be used to compute the
running error (estimates) of all nodes by traversing the
graph following the topological ordering, analogously to
the evaluation of the graph
in~\eqref{eq:fundop} and
the computation of the Jacobian
in~\eqref{eq:derprop}. This is implemented in the function \verb#eval_runerr#.

The first-order running error bound~\eqref{eq:run-err-bnd} is tighter than the
worst-case error bounds produced by an a priori first-order error analysis, but
can still be rather pessimistic for large graphs.
An alternative is to use the linearized version of
~\eqref{eq:run-err-lincomb},~\eqref{eq:run-err-mult},
and~\eqref{eq:run-err-ldiv} to estimate the value of $\delta_i$.
A stochastic estimate can be obtained, for instance, by interpreting $\eta$ as
a random variable uniformly distributed over the interval $[-u,u]$.

\begin{example}[Numerical cancellation]\label{exmp:roundoff}
  We illustrate the benefits of running error bounds on a classical example by
  Goldberg~\cite{gold91}: the evaluation of $y^2 - z^2$ for $y \approx z$, which is prone to cancellation if the two squares are not exactly representable in floating-point arithmetic.
  The following graph evaluates the scalar function$(1+2^{-27})^2-(1+2^{-28})^2$
  \vspace{-4pt}
\begin{lstlisting}
julia> x = 2. .^ -27
# Construct graph for (1+x)^2 - (1+x/2)^2.
julia> graph = Compgraph(Float64);
julia> add_lincomb!(graph, :y, 1, :I, 1, :x)
julia> add_lincomb!(graph, :z, 1, :I, 1/2, :x)
julia> add_mult!(graph, :y2, :y, :y)
julia> add_mult!(graph, :z2, :z, :z)
julia> add_lincomb!(graph, :out, 1, :y2, -1, :z2)
julia> add_output!(graph, :out)
\end{lstlisting}
  \vspace{-4pt}
The relative forward error and the corresponding running error bound are of
order $10^{-9}$ and $10^{-7}$, respectively. A more accurate probabilistic
running bound can also be computed.
\vspace{-4pt}
\begin{lstlisting}
# Compare forward error and running error bound.
julia> exact = 2. ^ (-2. * (i+1)) * (2. ^(i+2) + 3)
7.450580638557192e-9
julia> computed = eval_graph(graph, x, input = :x)
7.450580596923828e-9
julia> error = abs(exact - computed) / abs(exact)
5.5879354164678485e-9
julia> running_error_bound = eval_runerr(graph, x, input = :x)
2.98023227651711409552712846911887446840057366665064e-07
# A sharper stochastic bound can be found with the random mode.
julia> running_error_bound = eval_runerr(graph, x, input = :x, mode=:rand)
8.58378034267413297531189533682646018650457992705813e-09
\end{lstlisting}
\vspace{-4pt}
\end{example}

\section{Applications to state-of-the-art algorithms}\label{sec:use-cases}
We now show how the graph representation and the optimization described in previous sections can be used to improve the accuracy and performance of state-of-the-art algorithms.
As an illustration, we consider the computation of exponential and square root.

\paragraph{Matrix exponential}

The exponential is particularly well suited to our framework, as we can build on
recent developments in polynomial methods~\cite{sid19,sidr15,sast18,bbc19}.
In this context the optimization can be seen as the second step in a three-stage algorithm design procedure:
\begin{enumerate}[label=S\arabic*.,ref=S\arabic*]
\item\label{S1} Construct a graph using Paterson--Stockmeyer, monomials, or any method from~\cite{sid19,sidr15,sast18,bbc19}.
\item\label{S2} Optimize the coefficients in the graph with the technique in Section~\ref{sec:opt} run in high precision.
\item\label{S3} Compute guaranteed backward error bounds with the method discussed in Section~\ref{sec:bwerrexp}.
\end{enumerate}
For the function at hand, this approach brings an improvement over
state-of-the-art polynomial methods, which in several cases already
outperform  the Padé-based methods typically implemented in mathematical
software packages.

Implementations are named after the evaluation method used
(\verb#mono# for linear combination of monomials and
\verb#ps# for the Paterson--Stockmeyer scheme)
or after the authors of the paper where the algorithm is first presented:
\verb#sastre#
\cite{sast18},
\verb#sid#
\cite{sid19}, and
\verb#bbc#
\cite{bbc19}. More details on these implementations can be found in the online documentation of \verb#graph_sastre_exp# and \verb#graph_sid_exp#.
The suffix \verb#_opt# denotes improved graphs obtained by using our optimization strategy.
The same naming convention is used is the data repository, where the optimized graphs are available in CGR format together with the corresponding generated code for Julia, MATLAB, and C.

\newcommand{\mynum}[1]{#1}
\newcommand{\mynumsqr}[1]{\underline{#1}}
\newcommand{\mynumbf}[1]{\ensuremath{\mathbf{#1}}}

\begin{figure}[t]\centering
    \subfigure[$m=4$ multiplications.%
    \label{fig:sid_comparison_m4}]{%
      \pgfplotsset{
  every axis/.append style={title style={yshift=1ex}},
  grid=major,
  grid style={line width=.5pt, draw=Black!40},
  major grid style={dotted}}

\begin{tikzpicture}[trim axis left, trim axis right]

  \begin{axis}[
    xmin=0, xmax=0.8,
    ymin=1e-20, ymax=1e-15,
    xlabel={$x$},
    legend pos=south east,
    legend style={
      font=\tiny,
      draw=white!15!black,
      legend cell align={left},
      cells={line width=0.75pt}},
    ymode=log,
    width=2.8in,
    height=2.5in,
    cycle list name = bwd-list-3,
    every axis plot/.append style={thick}
    ]
    \addplot table[x={x}, y={sid}]{gfx/sid_comparison_m4.csv}; 
    \addplot table[x={x}, y={sid_opt}]{gfx/sid_comparison_m4.csv}; 
    \addplot coordinates {(0.69, 1e-30) (0.69, 1.1e-10)};
    \addplot coordinates {(0, 3.466574855304936e-17) (1, 3.466574855304936e-17)};
    \addplot coordinates {(0, 1.8e-17) (1, 1.8e-17)};
    \addlegendentry{\texttt{sid} (denoted by 15+ in \cite{sid19})};
    \addlegendentry{\texttt{sid\_opt}};
    \addlegendentry{Target disk radius $r=0.69$};
  \end{axis}
\end{tikzpicture}}%
    \qquad\qquad%
    \subfigure[$m=5$ multiplications.%
    \label{fig:sid_comparison_m5}]{%
      \pgfplotsset{
  every axis/.append style={title style={yshift=1ex}},
  grid=major,
  grid style={line width=.5pt, draw=Black!40},
  major grid style={dotted}}

\begin{tikzpicture}[trim axis left, trim axis right]

  \begin{axis}[
    xmin=0, xmax=2.5,
    ymin=1e-20, ymax=1e-15,
    xlabel={$x$},
    legend pos=south east,
    legend style={
      font=\tiny,
      draw=white!15!black,
      legend cell align={left},
      cells={line width=0.75pt}},
    ymode=log,
    width=2.8in,
    height=2.5in,
    cycle list name = bwd-list-3,
    every axis plot/.append style={thick}
    ]
    \addplot table[x={x}, y={sid}]{gfx/sid_comparison_m5.csv}; 
    \addplot table[x={x}, y={sid_opt}]{gfx/sid_comparison_m5.csv}; 
    \addplot coordinates {(1.9, 1e-30) (1.9, 1.1e-10)};
    \addplot coordinates {(0, 1.239149821803235e-16) (3, 1.239149821803235e-16)};
    \addplot coordinates {(0, 5.707080283800311e-17) (3, 5.707080283800311e-17)};
    \addlegendentry{\texttt{sid} (denoted by 21+ in \cite{sid19})};
    \addlegendentry{\texttt{sid\_opt}};
    \addlegendentry{Target disk radius $r=1.9$};
  \end{axis}
\end{tikzpicture}}%
    \caption{Relative error~\eqref{eq:relerr} over the positive real axis for the two methods \texttt{sid} and \texttt{sid\_opt} in Table~\ref{tbl:err}.
      The maximum error in $[0,r]$
      is marked with dotted lines.
      Clearly, the optimization is successful in the sense that it yields a smaller error in the domain $[0,r]$. As expected, for \texttt{sid\_opt} the error is more evenly distributed in contrast to \texttt{sid} which is very accurate near the origin. Since the coefficients of \texttt{sid} in \cite{sid19} were only given given to double precision, the coefficients in \texttt{sid\_opt} were also rouded to double precision for this figure.
      \label{fig:sid_comparison}
    }
\end{figure}
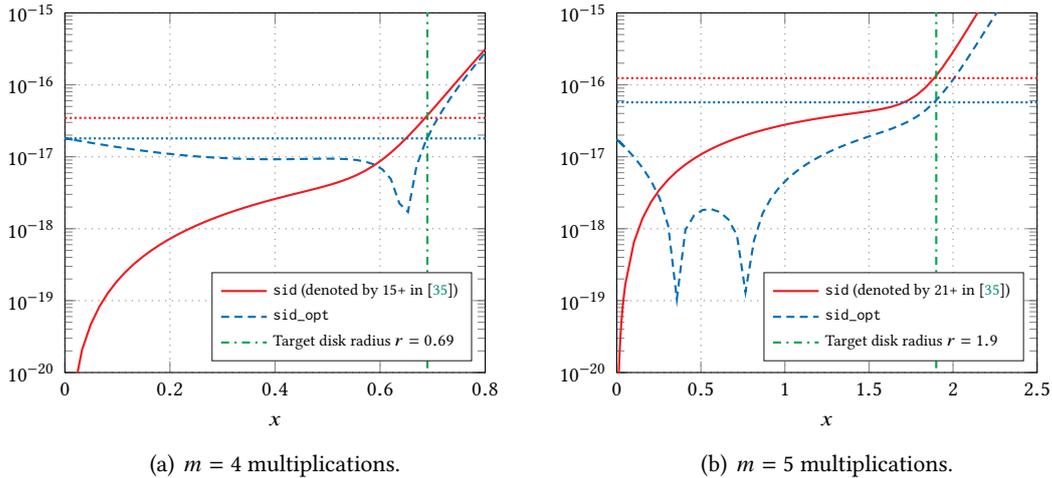

\begin{table}[t]
  \begin{center}
      \caption{Relative errors achieved after optimization within different disks with radius $r$ for $m$ multiplications.
       Optimization is carried out for the disk of radius $r$. The names of the graphs correspond to the graph files in the data repository.
      Underlined numbers correspond to optimizations initiated with squaring of the graph with one less multiplication.
      The dagger ${}^\dagger$ indicates that the graph in the corresponding columns were obtained using as starting graph for the optimization the graph in the previous column.
      \label{tbl:err}}
    \begin{tabular}{@{}
      l@{\extracolsep{4.3pt}}c@{\extracolsep{4.7pt}}
      c@{\extracolsep{4.7pt}}
      c@{\extracolsep{4.3pt}}c@{\extracolsep{4.7pt}}
      c@{\extracolsep{4.3pt}}c@{\extracolsep{4.7pt}}
      c@{\extracolsep{4.3pt}}c@{\extracolsep{4.7pt}}
      c@{}}
      \toprule
       &\multicolumn{1}{c}{$m=4$}
       &\multicolumn{2}{c}{$m=5$}
       &\multicolumn{2}{c}{$m=6$}
       &\multicolumn{2}{c}{$m=7$}
       &\multicolumn{1}{c}{$m=8$}
       \\
      \cmidrule(){2-2}
      \cmidrule(){3-4}
      \cmidrule(){5-6}
      \cmidrule(){7-8}
      \cmidrule(){9-9}
         \multicolumn{1}{c}{Graph} &
        $r=0.69$    &
       $r=1.68$  & $r=1.9{}^\dagger$  & $r=2.22$& $r=2.7{}^\dagger$& $r=3.59$& $r=6.4{}^\dagger$ & $\mynumsqr{r=13.5}$\\
      \midrule
      \multicolumn{9}{c}{Graphs found by optimization}\\
      \midrule
 \verb#mono_opt#  & $1.2 \!\cdot\!10^{-15}$ & $4.6 \!\cdot\!10^{-13}$ & $1.5 \!\cdot\!10^{-12}$ & $1.8 \!\cdot\!10^{-16}$ & $1.8 \!\cdot\!10^{-17}$ & $4.2 \!\cdot\!10^{-18}$ & \mynumbf{7.2 \!\cdot\!10^{-17}} & \mynumbf{3.7 \!\cdot\!10^{-17}}\\
 \verb#ps_opt#  & $1.2 \!\cdot\!10^{-15}$ & $1.8 \!\cdot\!10^{-13}$ & $2.0 \!\cdot\!10^{-12}$ & $5.6 \!\cdot\!10^{-14}$ & $3.3 \!\cdot\!10^{-15}$ & $2.5 \!\cdot\!10^{-20}$ & $3.3 \!\cdot\!10^{-13}$ & $5.0 \!\cdot\!10^{-14}$\\
 \verb#sastre_opt#  & $1.4 \!\cdot\!10^{-15}$ & $1.8 \!\cdot\!10^{-15}$ & \mynumbf{3.7 \!\cdot\!10^{-17}} & $1.5 \!\cdot\!10^{-16}$ & $6.0 \!\cdot\!10^{-17}$ & \mynumbf{1.3 \!\cdot\!10^{-21}} & $2.6 \!\cdot\!10^{-16}$ & $1.9 \!\cdot\!10^{-16}$\\
 \verb#bbc_opt#  & $2.3 \!\cdot\!10^{-14}$ & $4.8 \!\cdot\!10^{-14}$ & $1.3 \!\cdot\!10^{-13}$ & $\times$ & $\times$ & $\times$ & $\times$ & $\times$\\
 \verb#sid_opt#  & $1.1 \!\cdot\!10^{-16}$ & \mynumbf{1.8 \!\cdot\!10^{-16}} & $1.1 \!\cdot\!10^{-15}$ & \mynumbf{5.3 \!\cdot\!10^{-18}} & \mynumbf{1.6 \!\cdot\!10^{-18}} & $2.7 \!\cdot\!10^{-16}$ & $\times$ & $\times$\\
      \midrule
      \multicolumn{9}{c}{Polynomial methods from the literature} \\
      \midrule

      \verb#sastre#
                     &$1.5 \!\cdot\!10^{-11}$&  $\mynumsqr{5.8 \!\cdot\!10^{-14}}$ &$4.0 \!\cdot\!10^{-13}$&$1.8 \!\cdot\!10^{-8}$&$7.8 \!\cdot\!10^{-7}$&$\mynumsqr{6.4 \!\cdot\!10^{-10}}$&$\times$&$\times$\\
      \verb#bbc#
                     &$2.5 \!\cdot\!10^{-12}$&  $7.8 \!\cdot\!10^{-13}$&$9.9 \!\cdot\!10^{-12}$&$\times$&$\times$&$\times$&$\times$&$\times$\\
      \verb#sid#
       & $\mathbf{2.6 \!\cdot\!10^{-16}}$&$1.3 \!\cdot\!10^{-15}$&$5.1 \!\cdot\!10^{-15}$&$3.4 \!\cdot\!10^{-15}$&$6.0 \!\cdot\!10^{-14}$&$5.2 \!\cdot\!10^{-14}$&$5.4 \!\cdot\!10^{-8}$&$\times$\\
      \bottomrule
    \end{tabular}
  \end{center}
\end{table}
Table~\ref{tbl:err} and Figure~\ref{fig:sid_comparison} present the results of Stage \ref{S2}.
For various graphs, the table reports the scalar relative error
\begin{equation}\label{eq:relerr}
 \max_{z\in D(0,r)} \frac{|p(z)-\exp(z)|}{|\exp(z)|}=
\max_{z\in \partial D(0,r)} \frac{|p(z)-\exp(z)|}{|\exp(z)|},
\end{equation}
where $D(0,r)$ is the disk with radius $r$ centered at the origin,
and in practice we replace the boundary $\partial D(0,r)$ with a
sufficiently fine uniform discretization.
For comparison we have included the results for the optimized algorithms
as well as for some polynomial algorithms from the literature.
For the graphs found by optimization, the symbol $\times$ denotes combinations
for which the optimization routine gave unreasonably bad results; for the
graphs from the literature it denotes cases for which the coefficients of the
polynomial methods were not tabulated or directly computable. An example of the improvement
achieved by the optimization for two specific cases is also visible in Figure~\ref{fig:sid_comparison}.

In terms of relative error, the algorithm in \cite{sid19} can only be slightly improved for
$m=4$.
This is not surprising: the method is almost optimal, in the sense that it matches 15
Taylor coefficients out of the optimal 16 (cf. Table~\ref{tbl:deg}).
For $m>4$ multiplications, however, our approach brings a substantial improvement. Figure~\ref{fig:sid_comparison} shows that the methods that approximate the Taylor coefficients are the more accurate the closer $x$ is to the origin, whereas our approximation has a more balanced error on a disk.

All the graphs were represented in degree-optimal form
(Section~\ref{sec:degopt})
and all free coefficients in that form were used as free variables
in the optimization. For example,
the rows marked \verb#mono_opt#, \verb#ps_opt#, and \verb#sid_opt# correspond to
the degree-optimal forms in Figures~\ref{fig:monomial_degopt},~\ref{fig:ps_degopt},
and~\ref{fig:sid_degopt}, respectively.

\begin{remark}[Technical details for optimization]
  Even though for $m\le 5$ optimizing in double precision
  led to reasonably good results, in most cases we had to resort to
  higher precision. This was done using the Julia \verb#BigFloat# data type,
  which in turn relies on the GNU MPFR library~\cite{fhlp07}.
  Even when using higher precision, the complicated singular structure of
  the linear systems occurring during the optimization led to
  either stagnation or divergence, regardless of whether the backslash operator
  (based on the QR factorization) or the normal equations were used.
  We opted for a pseudoinverse (with adapted drop tolerance)
  computed using the high precision SVD routine in the Julia
  \verb#GenericLinearAlgebra.jl# package. The optimization routine requires
  a damping that has to be adapted throughout the iteration.
  We were unable to find a heuristic damping strategy
  suitable for all cases, thus we used scripting to decide and control
  damping, drop tolerances, and other parameters.
  In order to break the symmetry in the graph,
  we also sometimes added a random perturbation. The sequence of
  parameters that led
  to the results below are given in the companion data package.

  These problems stem from the singularity of the Jacobian, a consequence of
  the redundancy of the coefficients in the degree-optimal form.
  Removing obvious redundancies, for instance by fixing the first row of $H_a$
  and $H_b$, did not improve the result, and completely characterizing all
  redundancies in the degree-optimal form is a non-trivial issue
  beyond the scope of our package.
\end{remark}

\begin{table}[t]
  \caption{Values of $\thbwd$ in~\eqref{eq:thfwdi} for the guaranteed backward error bound computed with the method described in Section~\ref{sec:bwerrexp} for $m$ multiplications (or $m-1$ multiplications and one inverse multiplication). The values marked with ${}^*$ are those tabulated in the corresponding papers \cite{sid19} and \cite[Table~2]{bbc19}. A large value of $\thbwd$ indicates that the method is applicable in a wider domain. Boldface marks the best method for a given number of multiplications.
    \label{tbl:theta}
  }
  \begin{center}
    \begin{tabular}{@{}
      l@{\extracolsep{9.5pt}}c@{\extracolsep{14.5pt}}
      c@{\extracolsep{9.5pt}}c@{\extracolsep{14.5pt}}
      c@{\extracolsep{9.5pt}}c@{\extracolsep{14.5pt}}
      c@{\extracolsep{9.5pt}}c@{\extracolsep{14.5pt}}
      c@{}}
      \toprule
       &\multicolumn{1}{c}{$m=4$}
       &\multicolumn{2}{c}{$m=5$}
       &\multicolumn{2}{c}{$m=6$}
       &\multicolumn{2}{c}{$m=7$}
       &\multicolumn{1}{c}{$m=8$}
      \\
      \cmidrule(){2-2}
      \cmidrule(){3-4}
      \cmidrule(){5-6}
      \cmidrule(){7-8}
      \cmidrule(){9-9}
        \multicolumn{1}{c}{Graph}    &
        $r=0.69$    &
       $r=1.68$  & $r=1.9$  & $r=2.22$& $r=2.7$& $r=3.59$& $r=6.4$ & $13.5$\\
      \midrule
      \multicolumn{9}{c}{Graphs found from optimization} \\
      \midrule
 \verb#mono_opt# & 0.570 & 1.043 & 1.100 & 2.259 & 3.033 & 4.209 & \mynumbf{6.865} & \mynumbf{14.707}  \\
 \verb#ps_opt# & 0.569 & 1.142 & 1.134 & 1.713 & 2.431 & 5.064 & 5.210 & 10.768  \\
 \verb#sastre_opt# & 0.565 & 1.478 & \mynumbf{2.054} & 2.247 & 2.868 & 5.832 & 6.572 & 13.884  \\
 \verb#bbc_opt# & 0.450 & 1.256 & 1.377 & $\times$ & $\times$ & $\times$ & $\times$ & $\times$  \\
 \verb#sid_opt# & 0.674 & 1.683 & 1.766 & 2.581 & \mynumbf{3.313} & 3.633 & $\times$ & $\times$  \\
      \midrule
      \multicolumn{9
      }{c}{Polynomial methods from the literature}\\
      \midrule
      \verb#sid#$^*$
       & \mynumbf{0.695} & 1.683 & 1.683& 2.219 &2.219  & 3.53\hphantom{0} & 3.53\hphantom{0}  & $\times$ \\
      \verb#bbc#$^*$
       &  0.299& 1.09\hphantom{0} & 1.09\hphantom{0} & 2.22\hphantom{0}& 2.22\hphantom{0}&$\times$ &$\times$ &$\times$  \\
      \verb#exp_native#$^*$
       & 0.25\hphantom{0} &0.95\hphantom{0} &0.95\hphantom{0}& 2.10\hphantom{0} &2.10\hphantom{0}& 5.4\hphantom{00} &5.4\hphantom{00}& 10.8\hphantom{00}\\
      \bottomrule
    \end{tabular}
  \end{center}
\end{table}

Table~\ref{tbl:theta} reports the results of step~\ref{S3}, that is,
the application of the technique in Section~\ref{sec:bwerrexp}
to the graphs in Table~\ref{tbl:err}.
More precisely, the values correspond to the radii of the disks
that produced a backward error smaller
than the unit round-off $u$. The algorithm
underlying the graph is suitable for any matrix with norm
smaller than the value in the corresponding cell of the table. Note that
we optimize over one disk, and obtain a backward error guarantee for
another disk which may be smaller or larger.
Once again the optimization brings substantial improvements over the
algorithms used as initial guesses.
The best method is mostly consistent with the best method in
Table~\ref{tbl:err}, with the only exception of the case $m=5$.
The values in this table are found as in
Figure~\ref{fig:thbwdi_b}.

The table also reports the values of $\thbwd$ in~\eqref{eq:thbwdi} for the
Julia implementation of the Padé-based scaling-and-squaring algorithm
(\verb#exp_native#), which is based on the algorithm by Higham~\cite{high05e}.
This method always requires one multiplication by the
inverse of a matrix. In the table we assume that matrix multiplication and
multiplication by the inverse have same computational cost,
to the disadvantage of the methods we propose.

\newcommand{\ffrac}{\ensuremath{\phantom{\tfrac13}}}
\begin{table}[t]
  \caption{CPU timings and normalized computational cost $C$ for matrices
    of size $n=3000$.
    Method that are conclusively faster are in boldface.
      \label{tbl:cputime}}
  \begin{center}
    \begin{tabular}{@{}l@{\extracolsep{11.1pt}}l@{\extracolsep{20pt}}
      c@{\extracolsep{11.1pt}}c@{\extracolsep{20pt}}
      c@{\extracolsep{11.1pt}}c@{\extracolsep{20pt}}
      c@{\extracolsep{11.1pt}}c@{}}
      \toprule
      \multicolumn{8}{c}{Out-of-the-box}\\
      \midrule
      && \multicolumn{2}{c}{$\|A\|=2.5$}
      &\multicolumn{2}{c}{$\|A\|=6.0$}
      &\multicolumn{2}{c}{$\|A\|=13.5$}  \\
      \cmidrule{3-4}
      \cmidrule{5-6}
      \cmidrule{7-8}
      \multicolumn{1}{c}{Implementation} & \multicolumn{1}{c}{Language}
      & CPU time& $C$ & CPU time & $C$ &CPU time & C\\
      \midrule 
      \verb#expm# & MATLAB/MKL & 4.82 & $7\tfrac13$ & 5.43& $8\tfrac13$ &6.03&$9\tfrac13$  \\
      \verb#exp# & Julia/OpeBLAS & 5.62 & $7\tfrac13$ & 6.35 & $8\tfrac13$& 6.94&$9\tfrac13$ \\
      \verb#exp# & Julia/MKL & 4.81  & $7\tfrac13$ & 5.44 & $8\tfrac13$&6.04 &$9\tfrac13$ \\
      \verb#expmpol# \cite{sid19} & MATLAB/MKL & \mynumbf{4.30}& 7\ffrac & \mynumbf{4.94}  & 8\ffrac &5.52& 9\ffrac \\
      \midrule
      \multicolumn{8}{c}{Generated code}\\
      \midrule
      &&\multicolumn{2}{c}{$\|A\|=2.5$}
      &  \multicolumn{2}{c}{$\|A\|=6.0$}
      &\multicolumn{2}{c}{$\|A\|=13.5$} \\
      \cmidrule{3-4}
      \cmidrule{5-6}
      \cmidrule{7-8}
      \multicolumn{1}{c}{Graph} & \multicolumn{1}{c}{Language}
      & CPU time& $C$ & CPU time & $C$ & CPU time & $C$ \\
      \midrule
      Our approach& MATLAB/MKL    &           4.41       & 6\ffrac & 5.24           &7\ffrac& 6.13 &8\ffrac\\
      Our approach& Julia/OpenBLAS&           4.71       & 6\ffrac & 5.63           &7\ffrac& 6.48 &8\ffrac\\
      Our approach& Julia/MKL     & \mynumbf{3.80}       & 6\ffrac & \mynumbf{4.51} &7\ffrac& \mynumbf{5.24} &8\ffrac\\
      Our approach& C/MKL         &           4.03       & 6\ffrac & 4.84           &7\ffrac& 5.66 &8\ffrac\\
      \midrule
      \verb#exp_native_jl#& MATLAB/MKL    & 4.69 & $7\tfrac13$ & 5.28&$8\tfrac13$&         5.86 &$9\tfrac13$\\
\verb#exp_native_jl#& Julia/OpenBLAS& 5.28& $7\tfrac13$ & 5.84&$8\tfrac13$&$6.40$&$9\tfrac13$\\
\verb#exp_native_jl#& Julia/MKL  & 4.30& $7\tfrac13$ & \mynumbf{4.85}&$8\tfrac13$&\mynumbf{5.37}&$9\tfrac13$\\
\verb#exp_native_jl#& C/MKL         & 4.50 & $7\tfrac13$ & {5.04}&$8\tfrac13$&5.60&$9\tfrac13$\\
      \bottomrule
    \end{tabular}
  \end{center}
\end{table}

The actucal performance is illustrated by the CPU timings in
Table~\ref{tbl:cputime}.
The label \emph{Our approach} refers to
the corresponding graphs in Table~\ref{tbl:theta}: for $\|A\|=2.5$ we use \verb#sid_opt#, and for $\|A\|=6.0$ and $\|A\|=13.5$
we use \verb#mono_opt#. The matrices used for this experiment can be found in the data repository.
For reference we also report the CPU timing of the built-in functions for
computing the matrix exponential in Julia (\verb#exp#) and MATLAB
(\verb#expm#), 
and of the reference implementation of \verb#sid#\footnote{\url{http://personales.upv.es/~jorsasma/software/expmpol.m}}
run without norm estimation.
Among these \emph{out-of-the-box}
methods, the polynomial method in \verb#sid# is the fastest.
The column $C$ reports the asymptotic computational cost normalized so that
one $n \times n$ matrix multiplication has cost $C=1$ and the solution
of a linear system with $n$ right-hand sides has cost $4/3$.

In order to gauge the benefit of code generation, we also included
\verb#exp_native_jl#, the code generated from a graph implementation of the
built-in Julia \verb#exp# function.
By comparing the third row in Table~\ref{tbl:cputime} with the last four,
we see that the code generation is competitive, being even faster
than the (equivalent) native Julia implementation.

\begin{remark}[Technical details for comparisons of CPU timings]
  The timings were obtained following the best practices for
  benchmarking dense matrix computations. The processor boost was disabled
  in order to improve reproducibility.
  The reported values are median values of several runs.
  The experiments were run on a machine equipped with an 8-core Intel(R) Core(TM) i7-8550U CPU running at 1.80GHz. We used the development version 1.7.0 of Julia, as improvements made by the authors to the built-in Julia \verb#exp# function are only available in the development branch at the time of writing.
  The MKL implementation (version 2020.0.166) of the BLAS was used for the simulations in C.

  We used the largest matrices that would not cause the machine to run out of memory.
  The approach based on the degree-optimal form requires more memory than the other methods, as the number of intermediate matrices that have to be stored grows linearly with the number of multiplications.
  We also remark that our graph-based methods become more and more competitive as the matrices grow larger, in view of their lower computational cost $C$.
  In some languages we could make the matrices larger. In Julia with MKL, for instance, we could run the experiment for $n=4000$, obtaining for $\|A\|=13.5$ a timing of for our approach of 11.8s, compared with 13.6s for \verb#exp#, $13.1$s for \verb#exp_native_jl#, and 12.4s for  \verb#expmpol# in MATLAB.
\end{remark}

The results in this section are of empirical nature.
The optimization problem has many local minima, and a slight perturbation of the
initial parameters can lead the optimization to converge to a different point
of minimum.
There are several somewhat unexpected consequences of this. In the third row for $m=5$ in Table~\ref{tbl:err}, optimizing over a larger disk yields a smaller relative error. We were lucky. In Table~\ref{tbl:theta}, $m=5$ we observe
that \verb#sastre_opt#
has a larger $\theta$-value than \verb#sid_opt#, although the warmstart for the latter
corresponds to a more accurate approximation of the exponential.
These experimental results are of scientific value since both the graphs
and coefficients are publicly available.

\paragraph{Square root}
The approach, both in terms of optimization and
graph representation, can be applied to any analytic matrix function.
We illustrate this with the function $f(z)=\sqrt{z+1}$
in the domain $D(0,\frac12)$,  i.e., the disk centered at 0 with radius $\frac12$.
As this function is equivalent to the square root $\sqrt{z}$ in $D(1,\frac12)$, i.e., the disk centered at $1$ with radius $\frac12$, we can use
the Denman--Beavers iteration \cite{debe76} to warm start the optimization.
\begin{lstlisting}
julia> f=s->sqrt(s+1)
julia> (graph,cref)=graph_denman_beavers(2); # Use two iterations
julia> rename_node!(graph,:A,:A_shift,cref); # Introduce a shift
julia> add_lincomb!(graph,:A_shift,1.0,:A,1.0,:I);
julia> abs(eval_graph(graph,0.4im)-f(0.4im))
2.6944276809715006e-8
\end{lstlisting}
Directly optimizing the Denman--Beavers coefficients did not bring
any improvement.
However, the Denman--Beavers iteration
can be written in a form analogous to the degree-optimal form \eqref{eq:BBC_A}, where
the outer multiplication is replaced by a left division:
\begin{subequations}\label{eq:degoptgen}
  \begin{align}
    B_1&=I,\\
    B_2&=A,\\
  B_3&=(x_{a,1,1}B_1+x_{a,1,2}B_2)^{-1}(x_{b,1,1}B_1+x_{b,1,2}B_2),   \\
  B_4&=(x_{a,2,1}B_1+x_{a,2,2}B_2+x_{a,2,3}B_3)^{-1}(x_{b,2,1}B_1+x_{b,2,2}B_2+x_{b,2,3}B_3),\\
  B_5&=(x_{a,3,1}B_1+x_{a,3,2}B_2+\cdots +x_{a,3,4}B_4)^{-1}(x_{b,3,1}B_1+x_{b,3,2}B_2+\cdots +x_{b,3,4}B_4),\\
  B_6&=(x_{a,4,1}B_1+x_{a,4,2}B_2+\cdots +x_{a,4,5}B_5)^{-1}(x_{b,4,1}B_1+x_{b,4,2}B_2+\cdots +x_{b,4,5}B_5),\\
  r(A)&=y_1B_1+y_2B_2+y_3B_3+\cdots+y_6B_6.
\end{align}
\end{subequations}
In the notation of the Denman--Beavers example in Figure~\ref{fig:db} we have,
if the shift is included,
$B_3=X_0^{-1}=(A+I)^{-1}$,
$B_4=Y_1^{-1}=\bigl(\frac12I+\frac12 B_3\bigr)^{-1}$,
$B_5=X_1^{-1}=(I+\frac12A)^{-1}$, and
$B_6=Y_2^{-1}=\bigl(\frac12 \bigl(\frac12 (I+B_3)+B_5\bigr)\bigr)^{-1}$,
and  we can express the coefficients in the
modified degree-optimal form in Figure~\ref{fig:degoptgen-a}.
This formulation provides more degrees of freedom, which can be used
to fit the function more accurately.
We can use our optimization scheme, which is designed to work on
arbitrary graphs, on these coefficients in order to obtain a
better approximation. This produces the improved approximation
in Figure~\ref{fig:degoptgen-b},
where the largest coefficients are normalized so that all elements
in $H_a$ and $H_b$ have absolute value below 1.

The error of the approximation is shown in Figure~\ref{fig:sqrt_plus_one_err}.
We note that the optimized coefficients lead to an approximation error
of the order of the unit round-off $u$ in the target domain.
For comparison we also report the approximation error of the Paterson--Stockmeyer
evaluation of the degree-5 Taylor approximant as well as the degree-optimal
polynomial with coefficients optimized using
the Paterson--Stockmeyer method to warm start the optimization.
The scheme \eqref{eq:degoptgen} clearly provides a much more accurate
approximation. For fairness, the number of multiplications and multiplications
by an inverse are chosen so that the methods have the same asymptotic
computational cost $C$: the rational approximations in
Figure~\ref{fig:sqrt_plus_one_err} requires $4$ multiplications by an inverse,
whereas the polynomial approximations require $5$ matrix multiplications.

\begin{figure}[t]
  \centering
\subfigure[Coefficients for Denman--Beavers with shift.]{\label{fig:degoptgen-a}\scalebox{0.825}{
    \begin{minipage}{0.6\textwidth}
      \begin{align*}
        \left[\begin{array}{c|c}
    H_a & H_b
        \end{array}\right]&=
\left[\begin{array}{ccccc|ccccc}
  1    &1  &&&&   1&0   &    &   &       \\
\frac12&0   &\frac12&&&   1&0   &0    &   &       \\
1&\frac12&0 &0  &  & 1&0   &0    &0   &       \\
\frac14&0&\frac14&0&\frac12 & 1&0   &0    &0   &0       \\
  \end{array}\right]\\[6pt]
  y&=\begin{bmatrix}
\frac14&
\frac18&
0&
\frac14&
0&
\frac12
\end{bmatrix}
      \end{align*}\vspace{3pt}
    \end{minipage}
    }}
  \subfigure[Optimized coefficients.]{\scalebox{0.825}{\label{fig:degoptgen-b}
    \begin{minipage}{1.2\textwidth}
\begin{align*}
H_a&=
\begin{bmatrix}
1.0\phantom{0000000000000000} & \phantom{-}0.7610413081657074\hphantom{0} & \hphantom{-}0.0\hphantom{0000000000000000} & \hphantom{-}0.0\hphantom{00000000000000000} & 0.0 \\
1.0\phantom{0000000000000000} & -0.07127247848266706 & \phantom{-}0.9541397430184381\hphantom{0} & \hphantom{-}0.0\hphantom{00000000000000000} & 0.0 \\
1.0\phantom{0000000000000000} & \phantom{-}0.5152503668289159\hphantom{0} & -0.10150239642416546 & \hphantom{-}0.016950835467818528 & 0.0 \\
0.49421282325538585 & \phantom{-}0.05943417585870563 & \phantom{-}0.5017024742690791\hphantom{0} & \phantom{-}0.5276435956563207\hphantom{00} & 1.0
\end{bmatrix}\\[6pt]
H_b&=\begin{bmatrix}
1.0 & -0.011141376535527443 & \hphantom{-}0.0\phantom{000000000000000000} & \hphantom{-}0.0\hphantom{00000000000000000} & 0.0\hphantom{0000000000000000} \\
1.0 & -0.02575361014723645\phantom{0} & -0.004282806676350282\hphantom{0} & \hphantom{-}0.0\hphantom{00000000000000000} & 0.0\hphantom{0000000000000000} \\
1.0 & \phantom{-}0.03845623901890617\phantom{0} & \hphantom{-}0.026744368550352622\hphantom{0} & -0.029043605652931636 & 0.0\hphantom{0000000000000000} \\
1.0 & \phantom{-}0.12762852235945304\phantom{0} & -0.0015101375211018255 & -0.14486471812413823\hphantom{0} & 0.01583234295751802
\end{bmatrix}\\[6pt]
y^T&=\begin{bmatrix}
    \phantom{-}0.28457285753903816\hphantom{0} \\
    \phantom{-}0.11495283295141033\hphantom{0} \\
    -0.0344930890633602\hphantom{00} \\
    \phantom{-}0.34606871085737384\hphantom{0} \\
    -0.006442438999085179\\
    \phantom{-}1.4670072325444539\hphantom{00}
\end{bmatrix}
\end{align*}\vspace{3pt}
\end{minipage}}}

\caption{
  Coefficients for the scheme \eqref{eq:degoptgen} as approximation
  of $f(x)=\sqrt{x+1}$. \label{fig:degoptgen}
    }
\end{figure}

\colorlet{circlecolor}{Blue}

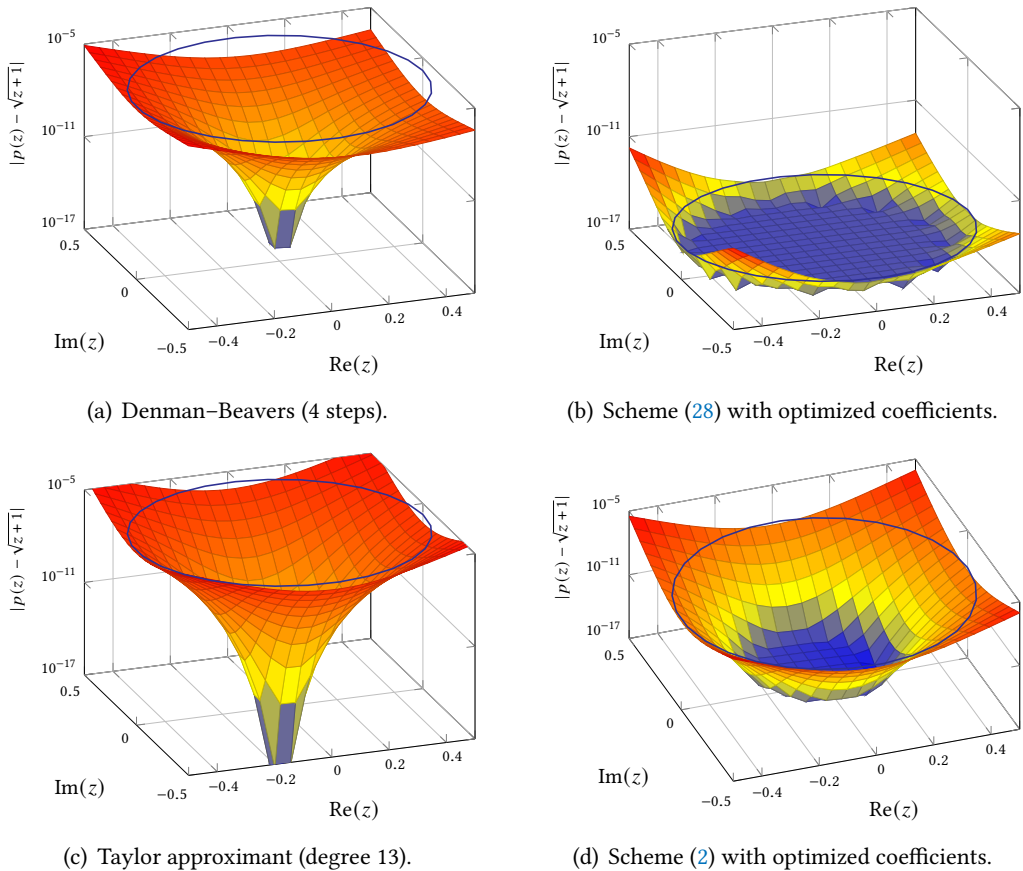
\begin{figure}[t]
  \begin{center}
    \subfigure[Denman--Beavers (4 steps).]{\scalebox{0.75}{
\tikzsetnextfilename{sqrt_plus_one_surf1}
\begin{tikzpicture}
\begin{axis}[view={-20}{30}, grid=both,zmode=log,
zmin=1e-17,zmax=1e-5,
xlabel={$\re(z)$},
ylabel={$\im(z)$},
zlabel={$|p(z)-\sqrt{z+1}|$}%
]
      \addplot3[surf] file {gfx/sqrt_plus_one_surf1.csv};
    \addplot3[color=circlecolor,thick] coordinates {
(0.500000,0.000000,  1.5e-6)
(0.488310,0.107485,  1.5e-6)
(0.453788,0.209945,  1.5e-6)
(0.398047,0.302587,  1.5e-6)
(0.323693,0.381081,  1.5e-6)
(0.234204,0.441756,  1.5e-6)
(0.133764,0.481775,  1.5e-6)
(0.027069,0.499267,  1.5e-6)
(-0.080891,0.493413, 1.5e-6)
(-0.185069,0.464488, 1.5e-6)
(-0.280594,0.413844, 1.5e-6)
(-0.362998,0.343850, 1.5e-6)
(-0.428429,0.257777, 1.5e-6)
(-0.473827,0.159651, 1.5e-6)
(-0.497069,0.054060, 1.5e-6)
(-0.497069,-0.054060,1.5e-6)
(-0.473827,-0.159651,1.5e-6)
(-0.428429,-0.257777,1.5e-6)
(-0.362998,-0.343850,1.5e-6)
(-0.280594,-0.413844,1.5e-6)
(-0.185069,-0.464488,1.5e-6)
(-0.080891,-0.493413,1.5e-6)
(0.027069,-0.499267, 1.5e-6)
(0.133764,-0.481775, 1.5e-6)
(0.234204,-0.441756, 1.5e-6)
(0.323693,-0.381081, 1.5e-6)
(0.398047,-0.302587, 1.5e-6)
(0.453788,-0.209945, 1.5e-6)
(0.488310,-0.107485, 1.5e-6)
(0.500000,-0.000000, 1.5e-6)
};
\end{axis}
\end{tikzpicture}
}}%
    \qquad%
    \subfigure[Scheme~\eqref{eq:degoptgen} with optimized
    coefficients.
    ]{\scalebox{0.75}{
\tikzsetnextfilename{sqrt_plus_one_surf2}
\begin{tikzpicture}
\begin{axis}[view={-20}{30}, grid=both,zmode=log,zmin=1e-17,zmax=1e-5,
xlabel={$\re(z)$},
ylabel={$\im(z)$},
zlabel={$|p(z)-\sqrt{z+1}|$}%
]
      \addplot3[surf] file {gfx/sqrt_plus_one_surf2.csv};
    \addplot3[color=circlecolor,thick] coordinates {
(0.500000,0.000000,  1.42e-15)
(0.488310,0.107485,  1.42e-15)
(0.453788,0.209945,  1.42e-15)
(0.398047,0.302587,  1.42e-15)
(0.323693,0.381081,  1.42e-15)
(0.234204,0.441756,  1.42e-15)
(0.133764,0.481775,  1.42e-15)
(0.027069,0.499267,  1.42e-15)
(-0.080891,0.493413, 1.42e-15)
(-0.185069,0.464488, 1.42e-15)
(-0.280594,0.413844, 1.42e-15)
(-0.362998,0.343850, 1.42e-15)
(-0.428429,0.257777, 1.42e-15)
(-0.473827,0.159651, 1.42e-15)
(-0.497069,0.054060, 1.42e-15)
(-0.497069,-0.054060,1.42e-15)
(-0.473827,-0.159651,1.42e-15)
(-0.428429,-0.257777,1.42e-15)
(-0.362998,-0.343850,1.42e-15)
(-0.280594,-0.413844,1.42e-15)
(-0.185069,-0.464488,1.42e-15)
(-0.080891,-0.493413,1.42e-15)
(0.027069,-0.499267, 1.42e-15)
(0.133764,-0.481775, 1.42e-15)
(0.234204,-0.441756, 1.42e-15)
(0.323693,-0.381081, 1.42e-15)
(0.398047,-0.302587, 1.42e-15)
(0.453788,-0.209945, 1.42e-15)
(0.488310,-0.107485, 1.42e-15)
(0.500000,-0.000000, 1.42e-15)
};
\end{axis}
\end{tikzpicture}
}}\\
    \subfigure[Taylor approximant (degree 13).]{\scalebox{0.75}{
\tikzsetnextfilename{sqrt_plus_one_surf3}
\begin{tikzpicture}
\begin{axis}[view={-20}{30}, grid=both,zmode=log,zmin=1e-17,zmax=1e-5,
xlabel={$\re(z)$},
ylabel={$\im(z)$},
zlabel={$|p(z)-\sqrt{z+1}|$}%
]
      \addplot3[surf] file {gfx/sqrt_plus_one_surf3.csv};
    \addplot3[color=circlecolor,thick] coordinates {
(0.500000,0.000000,  1.94e-6)
(0.488310,0.107485,  1.94e-6)
(0.453788,0.209945,  1.94e-6)
(0.398047,0.302587,  1.94e-6)
(0.323693,0.381081,  1.94e-6)
(0.234204,0.441756,  1.94e-6)
(0.133764,0.481775,  1.94e-6)
(0.027069,0.499267,  1.94e-6)
(-0.080891,0.493413, 1.94e-6)
(-0.185069,0.464488, 1.94e-6)
(-0.280594,0.413844, 1.94e-6)
(-0.362998,0.343850, 1.94e-6)
(-0.428429,0.257777, 1.94e-6)
(-0.473827,0.159651, 1.94e-6)
(-0.497069,0.054060, 1.94e-6)
(-0.497069,-0.054060,1.94e-6)
(-0.473827,-0.159651,1.94e-6)
(-0.428429,-0.257777,1.94e-6)
(-0.362998,-0.343850,1.94e-6)
(-0.280594,-0.413844,1.94e-6)
(-0.185069,-0.464488,1.94e-6)
(-0.080891,-0.493413,1.94e-6)
(0.027069,-0.499267, 1.94e-6)
(0.133764,-0.481775, 1.94e-6)
(0.234204,-0.441756, 1.94e-6)
(0.323693,-0.381081, 1.94e-6)
(0.398047,-0.302587, 1.94e-6)
(0.453788,-0.209945, 1.94e-6)
(0.488310,-0.107485, 1.94e-6)
(0.500000,-0.000000, 1.94e-6)
};
\end{axis}
\end{tikzpicture}
}}%
    \qquad
    \subfigure[Scheme~\eqref{eq:BBC_A} with optimized coefficients.]{\scalebox{0.75}{
\tikzsetnextfilename{sqrt_plus_one_surf4}
%
%
%
%
\begin{tikzpicture}
\begin{axis}[view={-20}{50}, grid=both,zmode=log,
         zmin=1e-17,zmax=1e-5,
xlabel={$\re(z)$},
ylabel={$\im(z)$},
zlabel={$|p(z)-\sqrt{z+1}|$}%
]
      \addplot3[surf] file {gfx/sqrt_plus_one_surf4.csv};
    \addplot3[color=circlecolor,thick] coordinates {
(0.500000,0.000000,  3.170000e-09)
(0.488310,0.107485,  3.170000e-09)
(0.453788,0.209945,  3.170000e-09)
(0.398047,0.302587,  3.170000e-09)
(0.323693,0.381081,  3.170000e-09)
(0.234204,0.441756,  3.170000e-09)
(0.133764,0.481775,  3.170000e-09)
(0.027069,0.499267,  3.170000e-09)
(-0.080891,0.493413, 3.170000e-09)
(-0.185069,0.464488, 3.170000e-09)
(-0.280594,0.413844, 3.170000e-09)
(-0.362998,0.343850, 3.170000e-09)
(-0.428429,0.257777, 3.170000e-09)
(-0.473827,0.159651, 3.170000e-09)
(-0.497069,0.054060, 3.170000e-09)
(-0.497069,-0.054060,3.170000e-09)
(-0.473827,-0.159651,3.170000e-09)
(-0.428429,-0.257777,3.170000e-09)
(-0.362998,-0.343850,3.170000e-09)
(-0.280594,-0.413844,3.170000e-09)
(-0.185069,-0.464488,3.170000e-09)
(-0.080891,-0.493413,3.170000e-09)
(0.027069,-0.499267, 3.170000e-09)
(0.133764,-0.481775, 3.170000e-09)
(0.234204,-0.441756, 3.170000e-09)
(0.323693,-0.381081, 3.170000e-09)
(0.398047,-0.302587, 3.170000e-09)
(0.453788,-0.209945, 3.170000e-09)
(0.488310,-0.107485, 3.170000e-09)
(0.500000,-0.000000, 3.170000e-09)
};
\end{axis}
\end{tikzpicture}
}}
    \caption{Error (forward) of four schemes for approximating the
      function $f(x)=\sqrt{x+1}$.
      The plots on the left correspond to 4 steps of the Denman--Beavers
      iteration (top) and to the Paterson--Stockmeyer evaluation of the
      truncated Taylor polynomial using $m=5$ multiplications (bottom).
      The plots on the right correspond to graphs with optimized coefficients;
      the computational graphs on the left in the same row were used as warm
      start for the optimization.
      The circles have radius 1/2 and $z$-value equal to the maximum
      relative error within the disk, which is
      $1.5\cdot 10^{-6}$,
      $1.4\cdot 10^{-15}$,
      $1.9\cdot 10^{-6}$ and
      $3.1\cdot 10^{-9}$
      for (a), (b), (c), and (d), respectively.
    \label{fig:sqrt_plus_one_err}
  }
  \vspace{-15pt}
\end{center}
\end{figure}

\section{Conclusions and outlook}\label{sec:conclusions}

We developed a Julia package for designing algorithms for matrix functions by working with computational graphs.
The software offers tools to generate and manipulate such graphs, to optimize their coefficients, and to produce code to evaluate them efficiently in several programming languages. Our numerical experiments confirm that the algorithm design technique we propose can improve state-of-the-art algorithms for evaluating the exponential and the square root of a matrix in many situations.

Our approach can be extended in several directions. The fixed-topology optimization is a general technique that can be applied to any matrix function, and can potentially be used to design computational methods for matrix functions for which efficient algorithms are currently not available. Moreover, the 3-step technique discussed in Section~\ref{sec:use-cases} can be extended to all those functions that can be evaluated efficiently using the scaling-and-squaring technique.

The degree-optimal polynomial form leads to a complete---though not minimal---parameterization of the space of polynomials that can be evaluated with a fixed number of non-scalar multiplications.
A minimal complete parameterization would be of interest as it would likely address the issues we experienced in the numerical simulations with singular Jacobians.
In general, establishing sufficient and practical conditions for a polynomial of degree at most $2^m$ to belong to the set of optimal-degree polynomial that can be evaluated with $m$ matrix multiplications remains an open problem.

In Section~\ref{sec:use-cases} we introduced a modified degree-optimal form \eqref{eq:degoptgen}, suitable for representing rational matrix functions, using which we obtained rational approximant to the square root much more efficient than polynomial approximant of equal computational cost. See \cite{sast12} for other approaches to evaluate rational approximants with few multiplications and left-divisions.
It may be worth investigating other variations of this already general form, such as, for example, sparse degree-optimal polynomials, where some coefficients are set to zero and not considered in the optimization.
These variants bring about a reduction of the memory footprint of the generated code.

By construction, algorithms obtained by using only~\ref{M1}--\ref{M3} evaluate rational functions. Our setting, however, is rather different from that of rational approximation theory, where the coefficients of a rational function of fixed degree are optimized.
The rational functions we consider have an additional constraint: a fixed topology or equivalently a \emph{fixed number of multiplications and divisions}, depending on the point of view.
Combining the effective methods for rational approximation with (say) the modified degree-optimal form in~\eqref{eq:degoptgen} could potentially lead to an improved workflow.

Although some tools for round-off error analysis are provided, further research would be required to understand which methods for the matrix exponential can guarantee a small round-off error.
The design and evaluation of graphs in Section~\ref{sec:use-cases} is based on high-precision arithmetic, but algorithms are typically run in binary32 or binary64 arithmetic. Therefore, despite being very rare in the literature, round-off error analysis of algorithms for matrix functions is of importance in practice.

\newcommand{\theacks}{We thank Jorge Sastre (Universitat Polit\`ecnica de Val\`encia) for interesting discussions and helpful comments and Olof Troeng (Lund universitet) for feedback on an early draft of the manuscript.}

\section*{Acknowledgments}\theacks

\bibliographystyle{plain-doi}
\bibliography{references}

\end{document}